\documentclass[12pt,twoside]{article}
\usepackage{amsmath}
\usepackage[all]{xy}
\expandafter\edef\csname :RestoreCatcodes\endcsname{%
   \catcode`\noexpand :=\the\catcode`:%
   \catcode`\noexpand @=\the\catcode`@%
   \catcode`\noexpand /=\the\catcode`/%
   \catcode`\noexpand &=\the\catcode`&%
   \catcode`\noexpand \^^M=\the\catcode`\^^M%
   \catcode`\noexpand \^^I=\the\catcode`\^^I%
   \expandafter\let\csname:RestoreCatcodes\endcsname=\noexpand\undefined}
\catcode`\:11  \catcode`\@11
   \let\:wlog\wlog \def\wlog#1{}
   \def\:wrn#1#2{\immediate\write\sixt@@n{--DraTeX warning--
      \ifcase #1
    DraTex.sty already loaded              
\or \string\Draw\space within \string\Draw 
\or Changing definition of \string#2      
\or No intersection points: #2            
\or Improper rotation of axes: #2         
\or (#2) in \string\DSeg\space is a point 
\fi}}
\def\:err#1#2{\errmessage{--DraTeX error-- \ifcase #1
     \string#2\space meaningless in three dimensions 
\or  \string#2\space meaningless in two dimensions   
\or  No \string\MarkLoc(#2)            
\or  \string#2 in three dimensions     
\or  Too many parameters in definition 
\or  \string\MoveFToOval(#2)? 
\fi}}
   \ifx\:Xunits\:undefined \else \:wrn0{} \fi
   \catcode`\ 9  \catcode`\^^M9 \catcode`\^^I9
      \newdimen\:LBorder \newdimen\:RBorder\chardef\:eight=8
\mathchardef\:cccvx=360
\newdimen\:mp    \:mp   0.1\p@
\newdimen\:mmp   \:mmp  0.01\p@
\newdimen\:mmmp  \:mmmp 0.001\p@
\newdimen\:XC    \:XC   90\p@
\newdimen\:CVXXX \:CVXXX180\p@
\newdimen\:CCCVX \:CCCVX\:cccvx\p@ \newdimen\:TeXLoc
\newbox\:box\newif\if:IIID  \newdimen\:Z   \newdimen\:Zunits
\newdimen\:Ex   \newdimen\:Ey  \newdimen\:Ez

\def\:AbsVal#1{ \ifdim#1<\z@-\fi #1 }
\def\:abs#1{\ifdim #1<\z@ #1-#1 \fi}
\def\:AbsDif#1#2#3{  #1#2   \advance#1  -#3
   \ifdim #1<\z@ #1-#1 \fi}
\def\:diff#1#2#3{ #1#2  \advance#1 -#3 }
\def\:average#1#2#3{
   #1#2  \advance#1  #3   \divide#1 \tw@}\def\:Opt#1#2#3#4{
   \def\:temp{
      \ifx      \:next\ifnum \def\:next{#3#1#4#2}
      \else\ifx \:next#1     \def\:next{#3}
      \else                  \def\:next{#3#1#4#2}\fi\fi \:next}
   \futurelet\:next\:temp}\def\Define#1{\:multid#1
   \:Opt(){\:Define#1}0}

\def\:DraCatCodes{\catcode`\ 9   \catcode`\^^M9
   \catcode`\^^I9  \catcode`\&13  \catcode`\~13 }

\def\:Define#1(#2){\begingroup  \:DraCatCodes  \::Define#1(#2)}

\def\::Define#1(#2)#3{\endgroup
   \let\:NextDefine\NextDefine
   \let\NextDefine\relax
   \ifcase#2\relax
      \def#1{#3}\or
      \:TxtPar\def#1(##1){#3}\or
         \:TxtPar\def#1(##1,##2){#3}\or
   \:TxtPar\def#1(##1,##2,##3){#3}\or
   \:TxtPar\def#1(##1,##2,##3,##4){#3}\or
   \:TxtPar\def#1(##1,##2,##3,##4,##5){#3}\or
   \:TxtPar\def#1(##1,##2,##3,##4,##5,##6){#3}\or
   \:TxtPar\def#1(##1,##2,##3,##4,##5,##6,##7){#3}\or
   \:TxtPar\def#1(##1,##2,##3,##4,##5,##6,##7,##8){#3}\or 
      \:TxtPar\def#1(##1,##2,##3,##4,##5,##6,##7,##8,##9){#3}\or
      \:err4{}\fi      \let\:TxtPar\relax  \:NextDefine}

\let\NextDefine\relax\let\:TxtPar\relax
\def\WarningOn{\def\:multid##1{
   \ifx ##1\:undefined \else \:wrn2##1\fi}}
\def\:gobble#1{}
\def\WarningOff{\let\:multid\:gobble}     \WarningOff
\Define\Indirect{\futurelet\:next\:Indirect}

\Define\:Indirect{\:theDoReg
   \ifx \:next<
      \def\:temp{\let\DoReg\:DoReg}
      \def\:next<##1>{\expandafter\:temp\csname :<##1>\endcsname}
   \else
      \def\:next##1<##2>{
         \expandafter\ifx \csname :<##2> \endcsname \relax
               \def\:next{##1}     \fi
\:indrwrn\Define     \:indrwrn\Object
\:indrwrn\Table      \:indrwrn\IntVar     \:indrwrn\DecVar 
         \def\:temp{\let\DoReg\:DoReg##1}
         \expandafter\:temp \csname :<##2> \endcsname}
   \fi      \:next}  \def\:indrwrn#1{    \def\:temp{#1}
   \ifx \:next\:temp \def\:wrn##1##2{\let\:wrn\::wrn} \fi}
\let\::wrn\:wrn 
\Define\:Hline{
   \setbox\:box\hbox{\vrule height0.5\:thickness
                    depth 0.5\:thickness width\:x}
   {\:d\:X \advance\:d \wd\:box
 \advance\:X -\:TeXLoc   \global\:TeXLoc\:d
 \vrule width\:X depth\z@ height\z@
 \raise \:Y \box\:box} }\Define\:Vline{
   \setbox\:box\hbox{\vrule width\:thickness
                 \ifdim \:y>\z@ height\:y  depth\z@
                 \else          height\z@  depth-\:y \fi}
   \advance\:X  -0.5\:thickness
   {\:d\:X \advance\:d \wd\:box
 \advance\:X -\:TeXLoc   \global\:TeXLoc\:d
 \vrule width\:X depth\z@ height\z@
 \raise \:Y \box\:box} }\Define\:MvTo(2){\:X#1\:Xunits \:Y#2\:Yunits}
\Define\:Mv(2){\advance\:X  #1\:Xunits
               \advance\:Y  #2\:Yunits}
\def\:DLn(#1,#2,{\:MvTo(#1,#2) \:LnTo(}\Define\:LnTo(2){
   \:x#1\:Xunits \advance\:x  -\:X
   \:y#2\:Yunits \advance\:y  -\:Y
   \:Ln(\:x\du,\:y\du) }\Define\:Ln(2){
   \:x#1\:Xunits \:y#2\:Yunits
   { \ifdim \:x<\z@
        \advance\:X \:x   \:x-\:x
        \advance\:Y \:y   \:y-\:y
     \fi
     \:yy\:AbsVal\:y
     \:dd\:yy \advance\:dd \:x
     \ifdim \:dd>\:mmmp
        \ifdim \:x>\:yy
           { \ifdim\:X<\:LBorder \global\:LBorder\:X\fi
\advance \:X  \:x
\ifdim \:X>\:RBorder \global\:RBorder\:X\fi }
           \let\:Yunitsy\:x  \let\:Xunitsx\:y  \let\:temp\:Hline
           \:dd0.6\:yy    \:divide\:dd\:x  
        \else
           { \advance \:X  -0.5\:thickness
\ifdim \:X<\:LBorder \global\:LBorder\:X \fi
\advance \:X  \:thickness
\advance \:X  \:x
\ifdim \:X>\:RBorder \global\:RBorder\:X\fi }
           \let\:Yunitsy\:y  \let\:Xunitsx\:x  \let\:temp\:Vline
           \:dd0.6\:x
\ifdim \:x>\:mmp    \:divide\:dd\:yy   \fi
        \fi
        \advance\:dd  0.4\p@
\:ragged\:Cons\:dd\:ragged
        \:HVLn
      \fi }
   \advance\:X  \:x   \advance\:Y  \:y}\Define\:HVLn{
   \:xx\:AbsVal\:Xunitsx  \:divide\:xx\:ragged
   \advance\:xx  0.99\p@   \:K\:InCons\:xx \relax
   \ifnum \:K>\z@
      \divide\:Xunitsx \:K  \advance\:K \@ne
      \divide\:Yunitsy \:K
   \else \:K\@ne  \fi  \:NextLn}

\Define\:NextLn{
   \ifnum\:K=\z@  \let\:NextLn\relax
   \else  { \:temp }  \advance\:K  \m@ne
      \advance\:X  \:x   \advance\:Y  \:y
   \fi   \:NextLn}\newdimen\:ragged

\Define\Ragged(1){    \:ragged#1\p@  \:ragged0.1\:ragged }
\Ragged(7.5)
\Define\PaintUnderCurve(4){{
   \:Z\:Y   \def\:next{\Curve(#1,#2,#3,#4)}
   \MoveToLoc(#1)  \:d\:X   \MoveToLoc(#4)
   \advance\:d -\:X
   \ifdim \:AbsVal\:d<\:mmp  \def\:next{}
   \else
      \def\:CrvLnTo(##1,##2){
         \:x \:X   \:y  \:Y    \:X\:DJ  \:Y\:yyyy
         \:xx\:X   \:dddd \:Y  \:X\:x  \:Y\:y
         { \advance\:Y  \:dddd   \divide\:Y \tw@
           \advance\:Z -\:Y
           \advance\:Y   0.5\:Z
           \:dddd  \:AbsVal \:Z  \:d\z@
           \def\:CrvLnTo{\:LnTo}
           \:yy\:Y  \:dd\:dddd  \:ddd\:dddd
           \::paint  }}
    \fi \:next       }}
\Define\DoCurve(1){ \let\:StopCurve\:SlowCurve
  \def\:CMv(##1){  \:x\:X \:y\:Y   \:MvTo(##1)
      \advance\:x -\:X   \advance\:y -\:Y
      \:xxx \:x    \:yyy\:y}
   \:DoCurve{\Curve(#1)}
   \let\:StopCurve\:FastCurve}

\def\:DoCurve#1(#2)#3{{\XSaveUnits
   \def\:next{#1}    \:MvTo(#2,#2)
   \:x\:AbsVal\:X  \:y\:Y  \:ddd\z@  \:length
   \:Z\:d   \:divide\:Z{1.41421\p@}
   \edef\:tempa{\the\:DoDist}   \global\:DoDist\z@
   \def\:CrvLnTo(##1){ \MarkLoc(1^)    \:CMv(##1)
      { \MarkLoc(2^)   \:ddd\z@    \:length
\:dd  \:DoDist  \global\advance\:DoDist  \:d
\:ddd \:DoDist  \:divide\:ddd\:Z
\DoReg\:InCons\:ddd  \:Z\DoReg\:Z

      \ifdim \:Z>\:dd
         \advance\:Z -\:DoDist
\advance\:dd -\:DoDist
\:divide\:Z\:dd
\advance\:X \:Cons\:Z\:xxx
\advance\:Y \:Cons\:Z\:yyy   \:DoRot 
         \def\:CrvLnTo{\:LnTo}
         \def\:OvalLn{\:Ln}  \XRecallUnits    #3 \fi}}
   \:next    \xdef\:DoDim{\:Cons\:DoDist}
   \global\:DoDist\:tempa     }
   \let\DoDim\:DoDim}

\newdimen\:DoDist\def\:DoRot{ \DSeg\RotateTo(1^,2^) }
\def\DoLine(#1,#2)(#3)#4{
   \MarkLoc($1)  \Move(#1,#2)
   \def\:next{   { \MarkLoc($2)
     \DSeg\RotateTo($1,$2)   \let\:DoRot\relax
     \edef\:RecallRagged{\the\:ragged} \MoveTo(#3,#3)
     \:x\:AbsVal\:X  \:y\:Y  \:ddd\z@  \:length
     \:ragged\:d   \divide\:ragged \tw@
     \DoCurve($1,$1,$2,$2)(#3)
        {\:ragged\:RecallRagged #4}  }
   \let\DoDim\:DoDim}  \:next }
\def\Table#1{\begingroup  \:DraCatCodes   \:multid#1
   \:DefineData#1}

\def\:DefineData#1#2{\endgroup
   \let\:temp~  \def~{\noexpand~}
   \edef#1{\noexpand\:DoPoly\expandafter\noexpand\csname :\string#1\endcsname}
   \expandafter\edef  \csname :\string#1\endcsname
       ##1{\noexpand\ifcase##1(#2)\noexpand\fi}
   \let~\:temp   \:DoNextPoly   \:DoNextPoly}

\def\:OR{\let\:or\or}  \:OR \catcode`\&13  \def&{)\noexpand\:or(}

\def\:TableData#1#2#3{\endgroup   \Table\:temp{#3}
   \:K\z@   \:J\z@    \def\:tempa(##1){\advance\:J \@ne }
   \:temp(0,999){\:tempa}    \let\:tempa&      \def\:temp{\def#1}
   \def&##1&{
      \ifnum  \:K<\:J
         \advance\:K \@ne
         \ifnum  \:K=\@ne   \def#1{#2(##1)}
         \else
            \:IIIexpandafter\:temp\expandafter{
               #1 & #2(##1) }
         \fi
      \else  \let\:next\relax \fi
      \:next}
   \let\:next&     &#3&&   \let&\:tempa }

\catcode`\&4  \def\:DoPoly#1(#2,#3)#4{
   \expandafter\let \csname :Back\the\:level\endcsname\:or
\expandafter\edef\csname :DoVars\the\:level\endcsname{
   \:DoB\the\:DoB}
\advance\:level  \@ne
   \:DoB#3  \advance\:DoB -#2
   \def\:PolyOr(##1){
      \ifnum  \:DoB=\z@  \:OR
      \else   #4(##1)   \advance\:DoB \m@ne    \fi}
   \:OR
   \def\:temp{\let\:or\:PolyOr #4}
   \:IIIexpandafter\:temp#1{#2}
   \advance\:level  \m@ne
\csname :DoVars\the\:level\endcsname
\def\:temp{\let\:or}
\expandafter\:temp\csname :Back\the\:level\endcsname }
\Define\PaintRect(2){\def\:next{{   \MarkLoc(^)
   \MoveToLoc(^) \Move( #1,0) \MarkLoc(^1) \Move(0,#2)
   \MarkLoc(^2)  \Move(-#1,0) \MarkLoc(^3)
   \PaintQuad(^,^1,^2,^3)        }}\:next}

\Define\PaintRectAt(4){\def\:next{{   \MoveTo(#1,#2)
   \:K#3   \advance\:K -#1
   \:DoB#4 \advance\:DoB -#2  \PaintRect(\:K,\:DoB)}}
   \:next}\def\::paint{
   \ifdim \:d<\:ragged      \advance\:xx -\:X
      \:yyy\:Y \:xxx\:dddd
      \advance\:Y \:yy  \divide\:Y \tw@
      \:average\:dddd\:dd\:ddd
      \def\:next{\:brush(\:xx,\z@)\:Y\:yyy\:dddd\:xxx}
   \else  \divide\:d \tw@
      \:average\:x\:X\:xx
      \:average\:y\:Y\:yy
      \:average\:dddd\:dd\:ddd
   \fi   \:next}\Define\:paint{{   \:AbsDif\:d\:xx\:x
   \ifdim \:d<\:mmp  \let\::paint\relax
   \else
      \ifdim  \:y >\:yyy   \:dd\:y   \:y \:yyy   \:yyy \:dd  \fi
      \ifdim  \:yy>\:yyyy  \:dd\:yy  \:yy\:yyyy  \:yyyy\:dd  \fi
      \:AbsDif\:dd\:yyyy\:yyy      \:AbsDif\:ddd\:yy\:y
      \ifdim \:dd<\:ddd  \:dd\:ddd  \fi
      \ifdim \:d >\:dd   \:d\:dd  \fi
      \advance\:d \:d      \:diff\:dd\:y\:yyy
      \:diff\:ddd\:yy\:yyyy
      \advance\:y    -0.5\:dd    \:abs\:dd
      \advance\:yy   -0.5\:ddd   \:abs\:ddd
      \:X\:x  \:Y\:y
      \:average\:dddd\:dd\:ddd
      \def\:next{ \:lpaint \:rpaint }
   \fi
   \::paint }}

\Define\:lpaint{ { \:xx\:x  \:yy \:y  \:ddd\:dddd  \::paint} }
\Define\:rpaint{ { \:X \:x  \:Y  \:y  \:dd \:dddd  \::paint} }
\Define\PaintQuad(4){\def\:next{{\Units(1pt,1pt)
   \MoveToLoc(#1)  \:x  \:X  \:y  \:Y
   \MoveToLoc(#2)  \:xx \:X  \:yy \:Y
   \MoveToLoc(#3)  \:xxx\:X  \:yyy\:Y
   \MoveToLoc(#4)  \:xxxx  \:X  \:yyyy \:Y
   
   \:paintQuad  }}\:next}

\def\:paintQuad{{
   \:SetVal\:a\:x\:y\:xx\:yy\:xxxx\:yyyy
\:SetVal\:b\:xx\:yy\:xxx\:yyy\:x\:y
\:SetVal\:c\:xxx\:yyy\:xx\:yy\:xxxx\:yyyy
\:SetVal\:cc\:xxxx\:yyyy\:xxx\:yyy\:x\:y
\def\:A{\:a} \def\:B{\:b} \def\:C{\:c} \def\:D{\:cc}
\:sort\:B\:A
\:sort\:C\:B  \:sort\:B\:A
\:sort\:D\:C  \:sort\:C\:B  \:sort\:B\:A
\let\:temp\relax
\:IsTriang\:A\:B
\:IsTriang\:B\:C
\:IsTriang\:C\:D
\:temp  
   \:Quad\:A\:temp>   \:xxxx\:xx \:yyyy\:yy \:Z\:d
\:Quad\:D\:next<   
   \:PrePaint(\:xx,\:yy,\:d,\:xxxx,\:yyyy,\:Z)
   \:temp  \:next }}
\Define\:PrePaint(6){
   \:x#1 \:y#2 \:yyy#3 \:xx#4 \:yy#5 \:yyyy#6
   \:paint }\def\:Quad#1#2#3{
         \:GetVal#1\:x\:y0
   \:GetVal#1\:xx\:yy1
   \:GetVal#1\:xxx\:yyy2
   \ifdim \:xx#3\:xxx
      \:ddd\:xx   \:xx\:xxx   \:xxx\:ddd
      \:ddd\:yy   \:yy\:yyy   \:yyy\:ddd
   \fi
          \def#2{}
   \:diff\:dd\:xxx\:xx
   \ifdim \ifdim\:AbsVal\:dd<\:mp  \:yy=\:yyy   \else  \z@>\z@  \fi
      
      \:d\:yy
   \else
      \ifdim \:AbsVal\:dd>\:mp
         \:diff\:dd\:xxx\:x   \:diff\:ddd\:yyy\:y
         \:divide\:ddd\:dd    \:diff\:dd\:xx\:xxx
         \:ddd\:Cons\:ddd\:dd   \advance\:ddd \:yyy
         \:d\:ddd
      \else \:d\:yyy \fi
      \edef#2{  \noexpand\:PrePaint
         (\the\:x ,\the\:y ,\the\:y,
         \the\:xx,\the\:yy,\the\:d)    }
   \fi}\def\:SetVal#1#2#3#4#5#6#7{
   \edef#1{(\the#2,\the#3,\the#4,\the#5,\the#6,\the#7)}}

\def\:sort#1#2{
   \ifdim \:IIIexpandafter\:field#1 <
          \:IIIexpandafter\:field#2
      \let\:temp#1  \let#1#2  \let#2\:temp
   \fi  }

\Define\:field(6){#1}

\def\:GetVal#1#2#3{
   \:IIIexpandafter\::GetVal #1#2#3}

\def\::GetVal(#1,#2,#3,#4,#5,#6)#7#8#9{
      \ifcase #9 #7#1   #8#2\or #7#3   #8#4\or #7#5   #8#6 \fi}
\def\:IsTriang#1#2{
   \ifdim \:IIIexpandafter\:field#1 =
          \:IIIexpandafter\:field#2
      \ifdim \:IIIexpandafter\:fieldB#1 =
             \:IIIexpandafter\:fieldB#2
         \def\:temp{ \:FixTria }
   \fi \fi  }

\def\:FixTria{
   \edef\:temp{\:IIIexpandafter\:FrsII\:B}
   \ifdim \:IIIexpandafter\:field\:A =
          \:IIIexpandafter\:field\:B
      \ifdim \:IIIexpandafter\:fieldB\:A =
             \:IIIexpandafter\:fieldB\:B
          \edef\:temp{\:IIIexpandafter\:FrsII\:C}
   \fi\fi
   \edef\:A{\:IIIexpandafter\:FrsII\:A}
   \edef\:D{\:IIIexpandafter\:FrsII\:D}
   \edef\:temp{
      \def\noexpand\:a{(\:A,\:temp,\:D)}
      \def\noexpand\:b{(\:temp,\:A,\:D)}
      \def\noexpand\:c{\noexpand\:b}
      \def\noexpand\:cc{(\:D,\:A,\:temp)}}
   \:temp
   \def\:A{\:a}  \def\:B{\:b}  \def\:C{\:c}  \def\:D{\:cc}  }

\Define\:fieldB(6){#2}
\Define\:FrsII(6){#1,#2}
\def\:IIIexpandafter{\expandafter\expandafter\expandafter}

\Define\DrawRect(2){ \Line( #1,0) \Line(0, #2)
                     \Line(-#1,0) \Line(0,-#2)}

\Define\DrawRectAt(4){{
   \MoveTo(#1,#2) \LineTo(#3,#2) \LineTo(#3,#4)
   \LineTo(#1,#4) \LineTo(#1,#2)}}
\def\du#1{ \ifx#1\:Xunits   \else\ifx#1\:Yunits
      \else\ifx#1\:Zunits   \else #1
      \fi \fi \fi}\Define\XSaveUnits{
   \expandafter\edef\csname XRecallUnits\the\:level\endcsname{
      \:StoreUnits}
    \advance\:level  \@ne}

\Define\XRecallUnits{
   \advance\:level \m@ne
   \csname XRecallUnits\the\:level \endcsname}

\Define\SaveUnits{
     }

\Define\:StoreUnits{       \:Xunits \the\:Xunits
   \:Yunits \the\:Yunits \:Zunits \the\:Zunits
   \:Xunitsx\the\:Xunitsx \:Xunitsy\the\:Xunitsy
\:Yunitsx\the\:Yunitsx \:Yunitsy\the\:Yunitsy }
\Define\:SearchDir{
   \ifdim \:x<\z@  \:x-\:x  \:y-\:y
   \edef\:tempA{\advance\:ddd  \ifdim \:y<\z@ - \fi\:CVXXX}
\else         \def\:tempA{} \fi
\ifdim \:y<\z@
   \edef\:tempA{\:ddd-\:ddd \advance\:ddd  \:CCCVX \:tempA}
   \:y-\:y
\fi
\ifdim \:y>\:x
   \:ddd\:y \:y\:x \:x\:ddd
   \edef\:tempA{\advance\:ddd  -\:XC \:ddd-\:ddd \:tempA}
\fi
   \:divide\:y\:x  \:d57.29578\:y
   \:ddd\:d       \:sqr\:y  \:K  \@ne
   \Do(1,30){\advance\:K \tw@
      \:d-\:Cons\:y\:d  \:dd\:d  \divide\:dd \:K
      \advance\:ddd \:dd  }
   \advance\:ddd -0.49182\:dd
   \:tempA }\Define\Curve(4){\def\:next{{  \XSaveUnits \Units(1pt,1pt)
   \MoveToLoc(#1) \:DI \:X  \:DK \:Y
   \MoveToLoc(#2) \:ddd \:X  \:Ez \:Y
   \MoveToLoc(#3) \:dd\:X  \:Ey\:Y
   \MoveToLoc(#4) \:DJ\:X  \:yyyy\:Y   \:Curve    }}\:next}
\Define\:FastCurve{   \:AbsDif\:d\:DI\:ddd  \:AbsDif\:dddd\:DK\:Ez
   \advance\:d \:dddd
   \ifnum \:d<\:ragged
      \:AbsDif\:d\:DJ\:dd \:AbsDif\:dddd\:yyyy\:Ey
      \advance\:d \:dddd                     \fi}
\let\:StopCurve\:FastCurve
\Define\:SlowCurve{
   \:AbsDif\:d\:DI\:DJ  \:AbsDif\:dddd\:DK\:yyyy
   \advance\:d \:dddd                     }
\Define\:Curve{   \:StopCurve
   \ifnum \:d<\:ragged   \:X\:DI \:Y\:DK
                    \:CrvLnTo(\:Cons\:DJ,\:Cons\:yyyy)
   \def\:SubCurves{}  \fi
   \:SubCurves}

\def\:CrvLnTo{\:LnTo}\Define\:SubCurves{
   \:average\:yy\:DI\:ddd     \:average\:Ex\:DK\:Ez
   \:average\:ddd\:ddd\:dd    \:average\:Ez\:Ez\:Ey
   \:average\:dd\:dd\:DJ  \:average\:Ey\:Ey\:yyyy
   \:average\:Zunits\:yy\:ddd    \:average\:Vdirection\:Ex\:Ez
   \:average\:ddd\:ddd\:dd    \:average\:Ez\:Ez\:Ey
   \:average\:DL\:Zunits\:ddd  \:average\:xxxx\:Vdirection\:Ez
   { \:ddd \:yy  \:Ez \:Ex   \:dd\:Zunits
     \:Ey\:Vdirection \:DJ\:DL \:yyyy\:xxxx \:Curve }
   \:DI \:DL \:DK \:xxxx               \:Curve   }
\def\MoveToCurve[#1]{
   \Define\:BiSect(3){\MoveToLoc(##1)
      \CSeg[#1]\Move(##1,##2)
   \MarkLoc(##3) }\:MvToCrv}

\Define\:MvToCrv(4){  \:BiSect(#1,#2,:a) \:BiSect(#2,#3,:b)
   \:BiSect(#3,#4,:c) \:BiSect(:a,:b,:A) \:BiSect(:b,:c,:B)
   \:BiSect(:A,:B,:Q)}
\Define\:OvalDir(3){   \:CosSin{#3\p@}
   \:Zunits\p@   \:d#2\:Zunits   \:x\:Cons\:d\:x
   \:Zunits\p@   \:d#1\:Zunits   \:y\:Cons\:d\:y
   \:SearchDir   }\Define\DrawOval   (2){ \DrawOvalArc(#1,#2)(0,\:cccvx) }
\Define\PaintOval  (2){ \PaintOvalArc(#1,#2)(0,\:cccvx) }
\Define\DrawCircle (1){ \DrawOval(#1,#1) }
\Define\PaintCircle(1){ \PaintOval(#1,#1) }
\def\DrawOvalArc(#1,#2)(#3,#4){{
       \:xxxx#4\p@  \advance\:xxxx -#3\p@
       \ifdim \:xxxx=\z@ \else
   \let\:SinOne\:SinB  \:OvalDir(#1,#2,#3)  \:DJ\:ddd
\:OvalDir(#1,#2,#4)  \:diff\:DI\:ddd\:DJ
\ifdim\:DI<\z@ \advance\:DI  \:CCCVX \fi
\ifdim \:xxxx<\:CCCVX \else \:DI\:CCCVX \fi
\:InitOval(#1,#2)  \:CosSin\:DJ 
   \:xxxx\:x  \:yyyy\:y   \:xx\:X  \:yy\:Y
   \advance\:X \:Xx\:x  \advance\:X \:Yx\:y
   \advance\:Y \:Xy\:x  \advance\:Y \:Yy\:y
   \let\:Xunits\empty  \let\:Yunits\empty
   \Do(1,\:InCons\:DI){
      \:dd\:X  \:ddd\:Y  \:X\:xx  \:Y\:yy
      \:AdvOv\:xxx\:yyy\:xxxx\:yyyy
      \:X\:dd  \:Y\:ddd
      \advance\:xxx -\:X   \advance\:yyy -\:Y
      \:d\:AbsVal\:xxx
      \advance\:d  \:AbsVal\:yyy
      \ifdim \:d>\:ragged
         \:OvalLn(\:xxx,\:yyy)  \fi  }
   \:OvalDir(#1,#2,#4)    \:CosSin\:ddd
   \advance\:xx \:Xx\:x  \advance\:xx \:Yx\:y
   \advance\:yy \:Xy\:x  \advance\:yy \:Yy\:y
   \advance\:xx -\:X     \advance\:yy -\:Y
   \:OvalLn(\:xx,\:yy) \fi }}

\def\:OvalLn{\:Ln}\def\DoOvalArc(#1)(#2){   \:xx\:X  \:yy\:Y
   \def\:CMv(##1){  \:Mv(\:xxx,\:yyy)
      \:x\:xxx  \:y\:yyy}
   \:DoCurve{             \:X\:xx  \:Y\:yy
       \def\:DoRot{  \let\:Xunits\:XunitsReg
                     \let\:Yunits\:YunitsReg
                     \DSeg\RotateTo(1^,2^)     }
       \let\::OvalLn\:CrvLnTo
\Define\:OvalLn(2){ \:dd\:AbsVal####1
   \advance\:dd \:AbsVal####2 \:divide\:dd\:ragged
   \:J\:InCons\:dd  \advance\:J  \@ne
   \divide####1  \:J  \divide####2  \:J
   \Do(1,\:J){\::OvalLn(####1,####2)}}
       \DrawOvalArc(#1)(#2)}}
\def\NextTable{\begingroup  \:DraCatCodes \:NextTable}

\def\:NextTable#1{\endgroup
  \def\:DoNextPoly{#1\NextTable{}}}
\NextTable{}\def\:AdvOv#1#2#3#4{
   \:d\:CosOne#3 \advance\:d -\:SinOne#4
   #4\:CosOne#4    \advance#4     \:SinOne#3   #3\:d
   \divide#3 \:eight  \divide#4 \:eight
   #1\:X   #2\:Y
   \:d\:Xx#3  \advance\:d \:Yx#4  \advance#1  \:d
   \:d\:Xy#3  \advance\:d \:Yy#4  \advance#2  \:d  }

\def\:CosOne{7.99878}   \def\:SinB{0.13962}
\def\PaintOvalArc(#1,#2)(#3,#4){{ \ifdim #3\p@=#4\p@
                      \let\:next\relax  \else
   \:d\:AbsVal{#1\:Xunits} \advance\:d \:AbsVal{#2\:Yunits}
   \ifdim \:d<3\:ragged   \divide \:d \tw@   \PenSize(\:d)
\:Mv(-0.5\:d\du,0)  \:Ln(\:d\du,0)
   \else    \:InitOval(#1,#2)
      \MarkLoc(o$)   \RotateTo(#3) \MoveFToOval(#1,#2)
      \:Ex\:X  \:Ey\:Y   \edef\:FirstOvDir{\:Cons\:ddd\p@}
      \MoveToLoc(o$) \RotateTo(#4) \MoveFToOval(#1,#2)
      \:Ez \:X  \:Vdirection \:Y   \edef\:LastOvDir{\:Cons\:ddd\p@}
      \MoveToLoc(o$)
      \if:rotated
      \:Zunits\p@     \:Zunits#1\:Zunits
      \:xx\:Cons\:Zunits\:Xunitsx
      \:Zunits\p@     \:Zunits#2\:Zunits
      \:yy\:Cons\:Zunits\:Yunitsx
      \:x\:xx  \:y\:yy  \:ddd\z@  \:length
      \:ddd\:d  \:divide\:xx\:ddd   \:divide\:yy\:ddd
\else \:xx\p@  \:yy\z@   \fi
      \:AbsDif\:d{#3\p@}{#4\p@}
      \ifdim \:d>359\p@ \:Ez-\:Xx\:xx  \advance\:Ez -\:Yx\:yy
\:Vdirection-\:Xy\:xx  \advance\:Vdirection -\:Yy\:yy
\advance\:Ez \:X  \advance\:Vdirection \:Y
\:setpaint\:PaintOvOv<> 
      \else
         \:x\:xx \:y\:yy \:SearchDir
\:xxx\:FirstOvDir  \advance\:xxx -\:ddd
\ifdim \:xxx<\z@    \advance\:xxx \:CCCVX  \fi
\:yyy\:LastOvDir   \advance\:yyy -\:ddd
\ifdim \:yyy<\z@    \advance\:yyy \:CCCVX  \fi
\:J\z@
\ifdim \:xxx<\:yyy  \ifdim        \:yyy<\:CVXXX
            \:Pntovln\:FirstOvDir\:FirstOvDir\:LastOvDir
\else \ifdim  \:xxx>\:CVXXX
            \:Pntovln\:FirstOvDir\:FirstOvDir\:LastOvDir
\else \:yyy-\:yyy  \advance\:yyy \:CCCVX
      \ifdim  \:xxx<\:yyy
            \:setpaint\:PntLeftOvOv><  \:PntMovln\:FirstOvDir\:xx
      \else \:FxLx  \:setpaint\:PntLeftOvOv><  \:Pntmovln\:LastOvDir
\fi  \fi  \fi
\else               \ifdim        \:xxx<\:CVXXX
            { \:setpaint\:PaintOvOv<> }  \:FxLx  \:setpaint\:PntLeftOvOv><
\:Usrch  \:PaintMidOvLn\:LastOvDir\:FirstOvDir
\else \ifdim  \:yyy>\:CVXXX
            {  \:FxLx \:setpaint\:PaintOvOv<>  }  \:setpaint\:PntLeftOvOv><
 \:Dsrch  \:PaintMidOvLn\:LastOvDir\:FirstOvDir
\else \:xxx-\:xxx  \advance\:xxx \:CCCVX
      \ifdim  \:yyy<\:xxx
            \:setpaint\:PaintOvOv<>    \:PntMovln\:FirstOvDir\:xx
      \else
            \:FxLx  \:setpaint\:PaintOvOv<>   \:Pntmovln\:LastOvDir
\fi  \fi  \fi  \fi
   \fi \fi \fi}}\def\:setpaint#1#2#3{{\aftergroup#1
\:d\:Xx\:xx   \advance\:d \:Yx\:yy
\ifdim \:d<\z@  \aftergroup#3
\else           \aftergroup#2  \fi}}
\def\:FxLx{\:d\:Ex  \:Ex\:Ez  \:Ez\:d}

\def\:Pntovln#1{
   \let\:SinOne\:SinB     \:CosSin#1
   \:xxxx\:Ex  \:yyyy\:Ey  \:Z\z@
   \:PaintOvLn}

\def\:PntMovln{
   \:Dsrch  \:FxLx   \:xx\:ddd  \:PaintOvLn}

\def\:Pntmovln{
   \:Usrch  \:FxLx  \:xx\:ddd   \:PaintOvLn\:xx}
\def\:PaintMidOvLn#1#2{
   \:FxLx   \:xx\:ddd  \:PaintOvLn#1#2
   \:xxx\:Ez  \:yyy\:Vdirection
   \:ddd\:Xy\:x  \advance\:ddd \:Yy\:y  \advance\:ddd \:Y
   \:PaintSlice}

\def\:PntLeftOvOv{
   \:xx-\:xx  \:yy-\:yy  \:PaintOvOv}

\def\:Dsrch{      \def\:SinOne{-\:SinB}
   \:xx\:x  \:yy\:y   \:SearchDir
   \:x\:xx  \:y\:yy
   \:d\:Ey  \:Ey\:Vdirection  \:Vdirection\:d }

\def\:Usrch{
   \let\:SinOne\:SinB
   \:x\:xx  \:y\:yy   \:SearchDir
   \:x\:xx  \:y\:yy}\def\:PaintOvLn#1#2{
   \:diff\:dd\:Ey\:Vdirection  \:diff\:ddd\:Ex\:Ez
   \ifdim \:AbsVal\:ddd>\:mp
      \:divide\:dd\:ddd      \:d#2
      \advance\:d -#1
      \ifnum \:d<\z@  \advance\:d \:CCCVX  \fi
      \:DoB\:InCons\:d  \let\:next\:PntDo  \:next
   \fi}

\def\:PntDo{
   \ifnum\:DoB=\z@ \let\:next\relax
   \else
      \::AdvOv\:x\:y
      \ifdim \:d>\:ragged
         \:ddd\:xxx  \advance\:ddd -\:Ez
         \:ddd\:Cons\:dd\:ddd
         \advance\:ddd \:Vdirection  \:PaintSlice
      \fi
      \advance\:DoB \m@ne
   \fi  \:next}\Define\:PaintSlice{   \:AbsDif\:dddd\:yyy\:ddd
   \advance\:yyy \:ddd  \divide\:yyy \tw@
   { \advance\:dddd \:Z  \divide\:dddd \tw@
      \:X\:xxxx  \:Y\:yyy
      \:xx\:xxx  \advance\:xxx -\:xxxx
      \:d\z@
      \:yy\:Y  \:dd\:dddd  \:ddd\:dddd
      \::paint     }
   \:xxxx\:xxx  \:yyyy\:yyy  \:Z\:dddd    \:J\z@}
\def\::AdvOv#1#2{  \:AdvOv\:xxx\:yyy#1#2
   \:AbsDif\:d\:xxxx\:xxx       \advance\:J \@ne
   \ifnum \:J=\sixt@@n   \multiply\:d \@cclvi
   \fi }\def\:PaintOvOv#1{   \def\:hdir{#1}
   \:xxx\:Xx\:xx  \advance\:xxx \:Yx\:yy
   \:yyy\:Xy\:xx  \advance\:yyy \:Yy\:yy
   \advance\:xxx \:X   \advance\:yyy \:Y
   \:Z\z@  \:xxxx\:xxx  \:yyyy\:yyy
   \:x\:xx   \:y\:yy   \:DoB\z@   \:J\z@
   \let\:next\:scanOvOv  \:next }
\Define\:scanOvOv{  \advance\:DoB \@ne
   \advance\:J \@ne
   \:AbsDif\:d\:xxx\:Ez
   \ifdim \ifdim      \:xxx\:hdir\:Ez    -
          \else\ifdim \thr@@\:d<\:ragged -
          \else\ifnum \:DoB>358            -
          \fi \fi \fi                     \p@<\z@
      \:ddd\:yyyy  \advance\:ddd -0.5\:dddd
      \:yyy\:yyyy   \advance\:yyy   0.5\:dddd
      \:xxx\:Ez \:PaintSlice
      \let\:next\relax
   \else
     \:d\:xx \advance\:d -\:x
     \:xxx\:yy \advance\:xxx -\:y
     \ifdim -\:Xx\:d\:hdir\:Yx\:xxx
        \let\:SinOne\:SinB
        \::AdvOv\:xx\:yy
        \ifdim  \ifdim       \:d >\:ragged -
                \else \ifnum \:J>\sixt@@n -
                \fi\fi                      \p@<\z@
           \:J\z@  \def\:SinOne{-\:SinB}
           \:AdvOv\:dd\:ddd\:x\:y  \:PaintSlice  \fi
     \else
        \def\:SinOne{-\:SinB}
        \::AdvOv\:x\:y
        \ifdim  \ifdim       \:d >\:ragged -
                \else \ifnum \:J>\sixt@@n -
                \fi\fi                      \p@<\z@
           \:J\z@   \let\:SinOne\:SinB
           \:AdvOv\:dd\:ddd\:xx\:yy  \:PaintSlice  \fi
   \fi \fi
   \:next}\Define\SetBrush{\:Opt[]\:SetBrush{}}

\def\:SetBrush[#1](#2,#3)#4{    \def\:temp{#4}
   \ifx \:temp\empty
      \def\:brush{   \let\:Xunits\empty \let\:Yunits\empty
                     \:thickness\:dddd  \:Ln  }
   \else       \def\:BruShape{#4}
      \:dd#2\:Xunits      \:ddd#3\:Yunits
      \edef\:Grd{ \:dd\the\:dd  \:ddd\the\:ddd }
      \MarkLoc($$)  \def\:temp{#1}
\ifx \:temp\empty  \:X\z@  \:Y\z@
\else  \MoveTo(#1) \fi
\edef\:BrOrg{ \:x\the\:X  \:y\the\:Y }
\MoveToLoc($$)
      \def\:brush(##1,##2){ \::brush }  \fi  }

\SetBrush(,){}\def\::brush{{  \SetBrush(,){}
\let\:Xunits\:XunitsReg   \let\:Yunits\:YunitsReg
\advance\:Y -0.5\:dddd   \:yy\:Y
\advance\:yy \:dddd
\:BrOrg  \:Grd  \advance\:xx  \:X
\ifdim \:xx<\:X  \:d\:X \:X\:xx \:xx\:d  \fi
   \:GridPt\:X\:x\:dd
   \:GridPt\:Y\:y\:ddd            \:x\:X
   \:DoBrush       }}\Define\:DoBrush{
   \ifdim       \:Y>\:yy  \let\:DoBrush\relax
   \else \ifdim \:X>\:xx \advance\:Y \:ddd     \:X\:x
   \else { \:BruShape }  \advance\:X \:dd
   \fi  \fi  \:DoBrush  }\def\:GridPt#1#2#3{   \:xxxx#1
   \advance#1 -#2  \:divide#1#3
   #1\:InCons#1#3  \advance#1 #2
   \ifdim #1=\:xxxx
   \else  \ifdim \:xxxx>#2 \advance#1 #3 \fi \fi  }
\newcount\:IntId    \newcount\:IntCount
\newcount\:DecId    \newcount\:DecCount

\Define\:NewCount{\alloc@ 0\count \countdef \insc@unt }
\Define\:NewDimen{\alloc@ 1\dimen \dimendef \insc@unt }

\def\:NewVar#1#2#3#4#5{ \:multid#1
   \def\:temp{    \csname \string#4\the#4\endcsname\z@
      \edef#1{\noexpand#3  \csname \string#4\the#4\endcsname}}
   \def\:next{  \global#5#4   \expandafter
       #2   \csname \string#4\the#4\endcsname  \:temp }
   \advance#4  \@ne
   \ifnum #4 > #5 \else \def\:next{\:temp}  \fi   \:next}

\def\IntVar#1{\:NewVar#1\:NewCount\:IntOp\:IntId\:IntCount}
\def\DecVar#1{\:NewVar#1\:NewDimen\:DecOp\:DecId\:DecCount}

\DecVar\Q  \DecVar\R   \DecVar\T
\IntVar\I  \IntVar\J   \IntVar\K
\def\WriteVal#1{\immediate\write\sixt@@n{...\string#1=#1;}}

\newdimen\:X   \newdimen\:Y
\newdimen\:x   \newdimen\:y   \newdimen\:d
\newdimen\:xx  \newdimen\:yy  \newdimen\:dd
\newdimen\:xxx \newdimen\:yyy \newdimen\:ddd
\newdimen\:xxxx\newdimen\:yyyy\newdimen\:dddd
\newcount\:J  \newcount\:K
\newdimen\:DI   \newdimen\:DJ
\newdimen\:DK   \newdimen\:DL   \newtoks\:t
\def\:IntFromPt#1#2{
   \:d#2\relax
   \advance\:d  \ifdim\:d<-0.5\p@-\fi  0.5\p@
   #1\:d    \divide#1  65536\relax}
\def\:temp{\catcode`\p12  \catcode`\t12}
\def\:Cons{\catcode`\p11  \catcode`\t11}
\:temp  \def\:Frac#1pt{#1}
        \def\:rnd#1.#2pt{#1}  \:Cons

\def\:Cons#1{\expandafter\:Frac\the#1}
\def\:sqr#1{#1\expandafter\:Frac\the#1#1}
\def\:InCons#1{\expandafter\:rnd\the#1}\def\:Val#1{#1;}
\let\Val\:Val\def\:IntOp#1#2{\csname :Op#2\endcsname#1}
\let\:SvIntOp\:IntOp
\def\:PreIntOp{\let\:IntOp\empty
   \let\Val\empty}
\def\:PostIntOp{\let\:IntOp\:SvIntOp
   \let\Val\:Val}

\expandafter\def\csname :Op;\endcsname#1{ \the#1}
\expandafter\def\csname :Op=\endcsname#1#2;{
   \:PreIntOp#1#2\:PostIntOp}
\expandafter\def\csname :Op+\endcsname#1#2;{
   \:PreIntOp\advance #1  #2\:PostIntOp}
\expandafter\def\csname :Op-\endcsname#1#2;{
   \:PreIntOp\advance #1  -#2\:PostIntOp}
\expandafter\def\csname :Op/\endcsname#1#2;{
   \:PreIntOp\divide#1   #2\:PostIntOp}
\expandafter\def\csname :Op*\endcsname#1#2;{
   \:PreIntOp\multiply#1   #2\:PostIntOp}
\def\:DecOp#1#2{ \csname :xOp#2\endcsname#1}
\let\:SvDecOp\:DecOp
\def\:PreDecOp{\let\:IntOp\the \def\:DecOp{\:Cons}
   \let\Val\empty   \let\:du\empty}
\def\:PostDecOp{\let\:IntOp\:SvIntOp \let\Val\:Val
   \let\:DecOp\:SvDecOp  \let\:du\::du  }

\def\::du#1{\p@
   \ifx#1\p@ \let\:temp\relax
   \else     \def\:temp{\du{#1}}
   \fi\:temp}                    \:PostDecOp

\expandafter\def\csname :xOp;\endcsname#1{ \:Cons#1}
\expandafter\def\csname :Op[\endcsname#1#2];{
   \:PreDecOp \:dd#2\p@  \:IntFromPt#1\:dd
                                 \:PostDecOp  }
\expandafter\def\csname :xOp=\endcsname#1#2;{
   \:PreDecOp#1#2\p@\:PostDecOp               }
\expandafter\def\csname :xOp(\endcsname#1#2){
   \:PreDecOp#1#2\p@\:PostDecOp               }
\expandafter\def\csname :xOp+\endcsname#1#2;{
   \:PreDecOp\advance #1  #2\p@\:PostDecOp  }
\expandafter\def\csname :xOp-\endcsname#1#2;{
   \:PreDecOp\advance #1  -#2\p@\:PostDecOp }
\expandafter\def\csname :xOp*\endcsname#1#2;{
   \:PreDecOp#1 #2#1\:PostDecOp              }
\expandafter\def\csname :xOp/\endcsname#1#2;{
   \:PreDecOp  \:divide#1{#2\p@}  \:PostDecOp }
\let\IF\ifnum

\def\EqText(#1,#2){
   \z@=\z@ \fi  \def\:temp{#1}
                \def\:next{#2}    \ifx \:temp\:next }

\def\:IfInt#1(#2,#3){ \z@=\z@ \fi
   \:IntOp\:K=#2;  \:IntOp\:J=#3; \ifnum  \:K#1\:J }

\def\:IfDim#1(#2,#3){ \z@=\z@ \fi
   \:DecOp\:d=#2;  \:DecOp\:dd=#3; \ifdim  \:d#1\:dd }

 \def\Do(#1,#2)#3{
   \expandafter\let
   \csname :Back\the\:level\endcsname\:Do
\expandafter\edef\csname :DoVars\the\:level\endcsname{
   \DoReg\the\DoReg \:DoB\the\:DoB}
\advance\:level  \@ne
   \DoReg#1  \:DoB#2  \relax
   \ifnum \DoReg<\:DoB
      \def\:Do{\ifnum \DoReg>\:DoB
                  \let\:Do\relax
               \else  #3\advance\DoReg  \@ne \fi
               \:Do}
   \else
      \def\:Do{\ifnum \DoReg<\:DoB
                  \let\:Do\relax
               \else  #3\advance\DoReg \m@ne  \fi
               \:Do}
   \fi  \def\:nextdo{ \:Do \advance\:level  \m@ne
\csname :DoVars\the\:level\endcsname
\def\:temp{\let\:Do}
\expandafter\:temp\csname
   :Back\the\:level\endcsname  } \:nextdo}

\let\:Do\relax

\newcount\DoReg   \let\:DoReg\DoReg       \newcount\:DoB
 \newcount\:level\def\::divide#1{   \:DI\:DK   \:dddd\:DL
   \advance\:DI -\:Cons\:dddd#1
   \:IntFromPt\:J\:dddd  \advance\:dddd -\:J\p@
   \multiply\:J  \@M   \:IntFromPt\:K{\@M\:dddd}
   \advance\:J \:K     \:dddd\@M\p@
   \divide\:dddd \:J   \advance#1 \:Cons\:DI\:dddd  }

\def\:divide#1#2{   \:DK#1   \:DL#2   #1\z@
   \::divide#1  \::divide#1  \::divide#1
   \::divide#1  \::divide#1  }
\def\:Sqrt#1{ \ifdim #1<\:mmp   #1\z@  \else
   \:dd#1   \divide\:dd \tw@
   \def\::Sqrt{  \:ddd#1
      \:divide\:ddd\:dd      \:AbsDif\:d\:dd\:ddd
      \advance\:dd \:ddd   \divide\:dd \tw@
      \ifdim  \:d < \:mmmp
         \let\::Sqrt\relax  \fi
      \::Sqrt}
   \::Sqrt   #1\:dd   \fi }\Define\:length{
   \:dd \:AbsVal \:ddd
   \:abs\:x  \:abs\:y
   \ifdim \:dd<\:x \:dd\:x \fi
   \ifdim \:dd<\:y \:dd\:y \fi
   \ifdim \:dd>\:mp
      \:divide\:x\:dd     \:sqr\:x
      \:divide\:y\:dd     \:sqr\:y
      \:divide\:ddd\:dd   \:sqr\:ddd
      \advance\:x  \:y  \advance\:x \:ddd
      \:y\:dd  \:Sqrt\:x  \:d\:Cons\:x\:y
   \else \:d\:dd \fi}\Define\:distance(2){\MarkLoc(@^)
   \MoveToLoc(#1)  \:x \:X  \:y \:Y  \:ddd\:Z
   \MoveToLoc(#2)
   \advance\:x  -\:X   \advance\:y  -\:Y
   \advance\:ddd  -\if:IIID \:Z \else \:ddd \fi
   \:length  \MoveToLoc(@^)}\def\:NormalizeDeg#1{
   \:DL#1   \:K\:InCons\:DL
   \divide\:K  \:cccvx   \multiply\:K  \:cccvx
   \advance #1 -\:K\p@
   \ifdim #1<\z@ \advance #1  \:CCCVX \fi
   \ifdim #1=\z@
      \ifdim\:DL=\z@ \else
         \advance #1  \:CCCVX \fi \fi }\def\:CosSin#1{ \:DK#1
   \:NormalizeDeg\:DK \def\:tempA{}
\ifdim \:CVXXX<\:DK
   \def\:tempA{\:y-\:y}
   \advance\:DK -\:CCCVX  \:DK-\:DK   \fi
\ifdim \:XC<\:DK
   \edef\:tempA{\:x-\:x \:tempA}
   \advance\:DK -\:CVXXX  \:DK-\:DK   \fi
\ifdim 45\p@<\:DK
   \edef\:tempA{\:d\:x \:x\:y \:y\:d \:tempA}
   \advance\:DK -\:XC   \:DK-\:DK   \fi
   \:x\p@   \:y0.01745\:DK   \:d\:y   \:K\@ne
   \edef\:next{\advance\:K \@ne
      \:sqr\:d  \divide\:d \:K  \advance}
   \:next \:x -\:d   \:next \:y -\:d
   \:next \:x  \:d   \:next \:y  \:d
   \:next \:x -\:d   \:next \:y -\:d
   \:next \:x  \:d   \:next \:y  \:d
   \:tempA   }   \Define\:rInitOval(2){
   \:Zunits\p@   \:dd#1\:Zunits
   \:d\:Cons\:dd\:Xunitsx   \edef\:Xx{\:Cons\:d}
   \:d\:Cons\:dd\:Xunitsy   \edef\:Xy{\:Cons\:d}
   \:dd#2\:Zunits
   \:d\:Cons\:dd\:Yunitsx   \edef\:Yx{\:Cons\:d}
   \:d\:Cons\:dd\:Yunitsy   \edef\:Yy{\:Cons\:d}}
\Define\:xyInitOval(2){
   \:d#1\:Xunits   \edef\:Xx{\:Cons\:d} \def\:Xy{0}
   \:d#2\:Yunits   \edef\:Yy{\:Cons\:d} \def\:Yx{0} }
 \def\:FigSize#1#2#3{
   \:x\:LBorder   \:y\:RBorder   \:d\:TeXLoc
   {\Object\:temp{#3}
    \setbox\:box\hbox{ \:temp
       \multiply\:x by \tw@  \multiply\:y by \tw@
       \xdef\:FSize{ \noexpand#1=\:Cons\:x;
                       \noexpand#2=\:Cons\:y;}}}
   \global\:LBorder\:x   \global\:RBorder\:y
   \global\:TeXLoc \:d
   \:FSize}
\expandafter\let \csname 0:Ln \endcsname\:Ln
\expandafter\def\csname 1:Ln \endcsname{
      \advance\:x -\:X  \advance\:y -\:Y
      \csname 0:Ln \endcsname(\:x,\:y)  }

\newcount\:ClipLevel   \:ClipLevel\@ne

\Define\Clip{\futurelet\:next\:Clip}

\Define\:Clip{
   \ifx \:next[  \expandafter\:DefClipOut
        \else    \expandafter\:DefClip     \fi }

\def\:DefClipOut[#1]{ \:DefClip(#1) }
\Define\:DefClip(1){  \def\:temp{#1}
   \ifx \:temp\empty
      \:ClipLevel\@ne     \def\:next{\let\:Ln}
      \expandafter\:next\csname 0:Ln \endcsname
   \else  \def\:temp{\::DefClip(#1)}  \fi  \:temp }
\Define\::DefClip(1){  \MarkLoc(^)
   \:x\:X  \:y \:Y  \Move(#1)
   \ifdim\:x>\:X  \:dd\:X  \:X\:x  \:x\:dd  \fi
   \ifdim\:y>\:Y  \:dd\:Y  \:Y\:y  \:y\:dd  \fi
   \advance\:ClipLevel  \@ne
   \expandafter\edef\csname \the
      \:ClipLevel :Ln \endcsname{
          \:xxx \the\:x   \:yyy \the\:y
          \:xxxx\the\:X  \:yyyy\the\:Y
          \ifx \:next[   \noexpand\:ClipOut
          \else          \noexpand\:ClipIn     \fi }
   \let\:Ln\:ClipLn   \MoveToLoc(^)  }\def\:ClipLn(#1,#2){
   \:x#1\:Xunits \:y#2\:Yunits
   {  \let\:Xunits\empty  \let\:Yunits\empty
   \advance\:x \:X   \advance\:y \:Y
   \ifdim \:x<\:X  \:dd\:X \:X\:x \:x\:dd
                   \:dd\:Y \:Y\:y \:y\:dd  \fi
   \:diff\:dd\:X\:x   \:diff\:ddd\:Y\:y
\:Z \:AbsVal \:dd
\advance\:Z  \:AbsVal\:ddd
\ifdim \:Z>\sixt@@n\p@
   \divide\:dd   128
   \divide\:ddd  128  \fi
\:Z\:Cons\:y\:dd
\advance\:Z -\:Cons\:x\:ddd
\ifdim \:dd<\z@
  \:dd-\:dd  \:ddd-\:ddd  \:Z-\:Z
\fi
   \csname \the \:ClipLevel :Ln \endcsname      } 
   \advance\:X \:x   \advance\:Y \:y  }\Define\:ClipIn{
   \def\:next{\let\:next}
   \expandafter\:next\csname \the
      \:ClipLevel :Ln \endcsname
   \advance\:ClipLevel \m@ne
   { \:ClipLeft\:xxxx \:ClipDown\:yyy \:ClipUp\:yyyy \:next }
   { \:ClipRight\:xxx \:ClipDown\:yyy \:ClipUp\:yyyy \:next  }
   { \:ClipDown\:yyyy \:next }
     \:ClipUp\:yyy    \:next }

\Define\:ClipOut{
   \:ClipLeft\:xxx      \:ClipRight\:xxxx
   \:ClipUp  \:yyyy     \:ClipDown \:yyy
   \def\:next{\let\:next}
   \expandafter\:next
      \csname \the\:ClipLevel :Ln \endcsname
   \advance\:ClipLevel \m@ne      \:next  }\def\:ClipLeft#1{
   \ifdim       \:x<#1  \:KilledLine
   \else \ifdim \:X<#1  \:X#1
      \ifdim \:dd>\:mmmp
         \:Y\:Cons\:ddd\:X  \advance\:Y \:Z
         \:divide\:Y\:dd
   \fi \fi \fi     \:CondKilLn  }

\def\:ClipRight#1{
   \ifdim       \:X>#1  \:KilledLine
   \else \ifdim \:x>#1  \:x#1
      \ifdim \:dd>\:mmmp
         \:y\:Cons\:ddd\:x  \advance\:y \:Z
         \:divide\:y\:dd
   \fi \fi \fi    \:CondKilLn }

\Define\:CondKilLn{
   \:d\:x  \advance\:d -\:X
   \ifdim \:d<\z@  \:d-\:d \fi
   \ifdim \:y<\:Y  \advance\:d \:Y \advance\:d -\:y
   \else           \advance\:d \:y \advance\:d -\:Y  \fi
   \ifdim  \:d<\:mmp   \:KilledLine \fi
   \ifdim  \:thickness=\z@ \:KilledLine \fi}

\Define\:KilledLine{
   \let\:ClipLeft\:gobble  \let\:ClipRight\:gobble
   \let\:ClipUp  \:gobble  \let\:ClipDown \:gobble
   \let\:next\relax
   \expandafter\def\csname 1:Ln \endcsname{}}\def\:ClipUp#1{
   \:AbsDif\:d\:y\:Y
   \ifdim  \:d<\:ragged
      \advance\:y  0.5\:thickness
\advance\:Y -0.5\:thickness
\ifdim       \:Y>#1  \:KilledLine
\else \ifdim \:y>#1
   \:thickness#1 \advance\:thickness -\:Y
   \advance\:Y  0.5\:thickness  \:y\:Y
\else
   \advance\:Y  0.5\:thickness
   \advance\:y -0.5\:thickness
\fi  \fi
\:dd\p@   \:ddd\z@
\def\:temp{  \:Z\:Y }  \:temp
   \else \let\:temp\relax
          \ifdim \ifdim\:Y<\:y\:Y\else\:y\fi >#1  \:KilledLine
   \else  \ifdim  \::ClipUp#1\:X\:Y
   \else  \ifdim  \::ClipUp#1\:x\:y
   \fi \fi \fi \fi   \:CondKilLn  }

\def\::ClipUp#1#2#3{
#3>#1   #3#1
\ifdim \:AbsVal\:ddd>\:mmmp
   #2\:Cons\:dd#3  \advance#2 -\:Z
   \:divide#2\:ddd
\fi  }\def\:ClipDown#1{   \:Ex2#1
   \:Flip\:y  \:Flip\:Y  \:ClipUp#1
   \:Flip\:y  \:Flip\:Y  \:temp }

\def\:Flip#1{   #1-#1  \advance#1 \:Ex  }
\Define\:SetDrawWidth{
   \hsize\:RBorder      \advance\hsize -\:LBorder
   \leftskip -\:LBorder \rightskip\z@}
\newdimen\:thickness   \:thickness0.75\p@

\Define\PenSize(1){\:thickness#1\relax}
\Define\:Draw{   \ifvmode \noindent\hfil\fi
   \global\:TeXLoc\z@
   \vbox\bgroup
           \begingroup
   \def\EndDraw{
           \endgroup   \:SetDrawWidth
         \egroup}
   \:DraCatCodes      \parindent\z@    \everypar{}
\leftskip\z@     \rightskip\z@    \boxmaxdepth\maxdimen
\let\FigSize\:FigSize
\def\Draw{\:wrn1{}} \:CommonIID   \:InDraw }
\Define\:CommonIID{\def\RotateTo{\:RotateTo}

\def\MoveF{\:MvF}\def\LineToLoc{\:LnToLoc}}\newdimen\:Xunits   \:Xunits\p@
\newdimen\:Yunits   \:Yunits\p@
\Define\:InDraw{\:Opt()\::InDraw{\:Xunits,\:Yunits}}
\Define\::InDraw(2){   \:X\z@   \:Y\z@
   \:Xunits#1 \relax  \:Yunits#2  \:AdjRunits
   \global\:LBorder\z@    \global\:RBorder\z@
   \:loadIID   \leavevmode }

\Define\:Resize(2){
   \:Xunits#1\:Xunits  \relax
   \:Yunits#2\:Yunits  \:AdjRunits  }
\Define\:Units(2){
   \:Xunits#1 \relax   \:Yunits#2    \:AdjRunits   }
\Define\:loadIID{
   \def\LineAt{\:DLn}
   \def\LineTo{\:LnTo}
   \def\MoveTo{\:MvTo}
   \def\Line{\:Ln}
   \def\Move{\:Mv}
   \def\MoveF{\:MvF}
   \:rotatedfalse\def\:InitOval{\:xyInitOval}
   
   \def\Units{\:Units}}

\def\DrawOn{\def\Draw{\:Draw}}                      \DrawOn
\def\DrawOff{\def\Draw{\begingroup \:J\@cclv
                       \:NoDrawSpecials \:NoDraw}}
\catcode`\/0 \catcode`\\11
/def/:NoDraw#1\EndDraw{/endgroup}
/catcode`/\0 /catcode`//12

\def\:NoDrawSpecials{\catcode\:J11
  \ifnum \:J=\z@
     \let \:NoDrawSpecials\relax \fi
  \advance\:J  \m@ne \:NoDrawSpecials}
\let\:XunitsReg\:Xunits   \let\:YunitsReg\:Yunits  \Define\EntryExit(4){
   \edef\:InOut##1{
      \noexpand\ifcase ##1\space
         #1\noexpand\or #2\noexpand\or
         #3\noexpand\or #4\noexpand\fi}}

\EntryExit(0,0,0,0)\Define\:DrawBox{   \:x0.5\wd\:box
   \:y\ht\:box \advance\:y  \dp\:box
   \divide\:y \tw@
   \advance \:X -\:InOut0\:x
   \advance \:Y -\:InOut1\:y
   { \advance \:X -\:x   \advance \:Y -\:y
     { \ifdim\:X<\:LBorder \global\:LBorder\:X\fi
\advance\:X \wd\:box
\ifdim\:X>\:RBorder \global\:RBorder\:X\fi }
     \advance \:Y  \dp\:box
     {\:d\:X \advance\:d \wd\:box
 \advance\:X -\:TeXLoc   \global\:TeXLoc\:d
 \vrule width\:X depth\z@ height\z@
 \raise \:Y \box\:box} }
   \edef\MoveToExit(##1,##2){
   \:X\the\:X   \:Y\the\:Y
   \:x\the\:x   \:y\the\:y
   \advance\:X  ##1\:x
   \advance\:Y  ##2\:y}
   \advance \:X  \:InOut2\:x
   \advance \:Y  \:InOut3\:y }
 \Define\ThreeDim{\:Opt[]\:ThreeDim{\p@}}

\def\:ThreeDim[#1](#2){\::ThreeDim[#1](#2,,)}

\def\::ThreeDim[#1](#2,#3,#4,#5){ \bgroup\begingroup
   \def\EndThreeDim{          \endgroup\egroup}
   \:IIIDtrue   \:Zunits#1  \:Z\z@
   \def\:temp{#4}
   \ifx \:temp\empty   \:CosSin{#3\p@}  \:divide\:x\:y       \:Ey\:x
\:CosSin{#2\p@}  \:Ex\:Cons\:Ey\:x  \:Ey\:Cons\:Ey\:y
\let\:project\:projectPar

   \else               \:Ex#2\:Xunits \:Ey#3\:Yunits \:Ez#4\:Zunits
\let\:project\:projectPer
 \fi
   \def\LineAt{\:tDLn}
\def\LineTo{\:tLnTo}
\def\MoveTo{\:tMvTo}
\def\Line{\:tLn}
\def\Move{\:tMv}

\def\Units{\:tUnits}\def\RotateTo{\:tRotateTo}

\def\MoveF{\:tMvF}
\:Vdirection\z@ \def\LineToLoc{\:tLnToLoc} }  \Define\:projectPer{
   \:diff\:x\:X\:Ex   \:diff\:y\:Y\:Ey
   \:diff\:xxxx\:Z\:Ez
   \:divide\:x\:xxxx      \:divide\:y\:xxxx
   \:x-\:Cons\:Ez\:x     \:y-\:Cons\:Ez\:y
   \advance\:x  \:Ex    \advance\:y \:Ey}
\Define\:projectPar{
   \:x\:Cons\:Ex\:Z  \advance\:x \:X
   \:y\:Cons\:Ey\:Z  \advance\:y \:Y}
\def\:tDLn(#1,#2,#3,{\:tMvTo(#1,#2,#3)
                     \:tLnTo(}
\Define\:tMvTo(3){   \:X#1\:Xunits
   \:Y#2\:Yunits    \:Z#3\:Zunits}
\Define\:tMv(3){     \advance\:X  #1\:Xunits
   \advance\:Y  #2\:Yunits
   \advance\:Z  #3\:Zunits}
\Define\:tResize(3){   \:Xunits#1\:Xunits
   \:Yunits#2\:Yunits \:Zunits#3\:Zunits
   \:AdjRunits}
\Define\:tUnits(3){    \:Xunits#1  \relax
   \:Yunits#2  \relax       \:Zunits#3
   \:AdjRunits}\Define\:tLnTo(3){
   \:project   \edef\:temp{\:x\the\:x \:y\the\:y}
   \:X#1\:Xunits  \:Y#2\:Yunits  \:Z#3\:Zunits
   \:DLN}

\Define\:tLn(3){
   \:project   \edef\:temp{\:x\the\:x \:y\the\:y}
   \advance\:X  #1\:Xunits
   \advance\:Y  #2\:Yunits
   \advance\:Z  #3\:Zunits   \:DLN}

\Define\:DLN{
   \TwoDim      \:temp        \LineTo(\:x\du,\:y\du)
   \EndTwoDim }\Define\TwoDim{\bgroup\begingroup
   \def\EndTwoDim{\endgroup\egroup}
   \:loadIID
   \if:IIID  \:IIIDfalse \:project \:X\:x \:Y\:y
             \:CommonIID \fi
   \Units(\:Xunits,\:Yunits)}  
\newif\if:rotated
\Define\RotatedAxes(2){\begingroup
   
   \if:IIID  \:IIIDfalse
             \:project \:X\:x \:Y\:y
             \:CommonIID    \fi
   \:DK#2\p@ \advance\:DK -#1\p@
\advance\:DK -\:CVXXX   \:NormalizeDeg\:DK
\ifdim \:DK<\:mmp \:wrn4{(#1,#2)}
\else
  \:DK#1\p@ \:NormalizeDeg\:DK    \:CosSin\:DK
  \:Xunitsx\:x  \:Xunitsy\:y
  \:DK#2\p@ \advance\:DK -\:XC \:NormalizeDeg\:DK \:CosSin\:DK
  \edef\Units(##1,##2){\noexpand\:Units(##1,##2)
    \:Xunitsx \:Cons\:Xunitsx\:Xunits
    \:Xunitsy \:Cons\:Xunitsy\:Xunits
    \:Yunitsx-\:Cons\:y\:Yunits
    \:Yunitsy \:Cons\:x\:Yunits  } \fi
 \def\MoveTo{\:rMvTo}
\def\Move{\:rMv}
\def\LineTo{\:rLnTo}
\def\Line{\:rLn}
\def\MoveF{\:rMvF}\def\:InitOval{\:rInitOval}

   \MarkLoc(:org) \Units(\:Xunits,\:Yunits)  \:rotatedtrue }
\newdimen\:Xunitsx  \newdimen\:Xunitsy
\newdimen\:Yunitsx  \newdimen\:Yunitsy
\Define\:rResize(2){
   \:Xunits #1\:Xunits    \:Yunits #2\:Yunits
   \:Xunitsx#1\:Xunitsx   \:Xunitsy#1\:Xunitsy
   \:Yunitsx#2\:Yunitsx   \:Yunitsy#2\:Yunitsy  }
\Define\:rMvTo{\MoveToLoc(:org) \:rMv}
\Define\:rMv(1){ \:rxy(#1)
   \advance\:X  \:x \advance\:Y  \:y}
\Define\:rLnTo(1){ \MarkLoc(:a) \:rMvTo(#1) {\LineToLoc(:a)} }
\Define\:rLn(1){ \:rxy(#1) \:Ln(\:x\du,\:y\du) }
\Define\:rxy(2){ \:Zunits\p@   \:Zunits#1\:Zunits
   \:x\:Cons\:Zunits\:Xunitsx
   \:y\:Cons\:Zunits\:Xunitsy
   \:Zunits\p@   \:Zunits#2\:Zunits
   \advance\:x  \:Cons\:Zunits\:Yunitsx
   \advance\:y  \:Cons\:Zunits\:Yunitsy}
\Define\:rMvF(1){  \:CosSin\:direction
   \:Zunits\p@            \:Zunits#1\:Zunits
   \:x\:Cons\:Zunits\:x   \:y\:Cons\:Zunits\:y
   \edef\:temp{(\:Cons\:x,\:Cons\:y)}   \expandafter\Move\:temp}
\Define\:AdjRunits{
   \:Xunitsx\:Xunits   \:Xunitsy\z@
   \:Yunitsx\z@        \:Yunitsy\:Yunits}
\Define\Text{  \setbox\:box
   \vtop\bgroup    \edef\DoReg{\the\DoReg}
      \hyphenpenalty\@M  \exhyphenpenalty\@M
      \catcode`\ 10 \catcode`\^^M13 \catcode`\^^I10
      \catcode`\&4  \let~\space
      \:Text}                       \catcode`\^^M13 %
\def\:Text(--#1--){%
      \:SetLines#1\hbox{}^^M--)^^M %
   \egroup                  %
   \if:IIID \TwoDim  \:DrawBox  \EndTwoDim %
   \else             \:DrawBox  \fi}       %
\def\:SetLines#1^^M{        %
   \def\:TextLine{#1}       %
   \ifx \:TextLine\:LastLine   \let\:temp\relax      %
   \else  \def\:temp{                                 %
             \:IndirectLines#1\relax~~--)~~\:SetLines}%
   \fi  \:temp }                      \catcode`\^^M9
\def\:IndirectLines#1~~{    \def\:TextLine{#1}
  \ifx \:TextLine\:LastLine   \let\:temp\relax
  \else  \def\:temp{\:AddLine{#1}\:IndirectLines}
  \fi  \:temp }

\Define\:LastLine{--)}

\def\:AddLine#1{
   \ifvmode \noindent    \hsize\z@ \else
      \hfil \penalty-500 \hbox{}    \fi
   \hfil#1
   \setbox\:box\hbox{#1}
   \ifdim \wd\:box>\hsize \hsize\wd\:box \fi}\def\TextPar#1#2{
   \def\:TxtPar##1(##2){##1(--##2--)}
   \edef\:temp{\expandafter\noexpand\csname :\string#2\endcsname}
   \edef#2{\noexpand\:TextPar\expandafter\noexpand\:temp}
   \expandafter\let\:temp\:undefined
   \expandafter#1\:temp}
                                   \catcode`\^^M13
\def\:TextPar#1{\begingroup     \catcode`\&4         %
   \catcode`\ 10 \catcode`\^^M13 \catcode`\^^I10 %
   \:TPar{#1}}                     \catcode`\^^M9  %

\def\:TPar#1(--#2--){\endgroup
   #1(--#2--)  }
 \newdimen\:direction  \newdimen\:Vdirection
\Define\:Rotate(1){   \advance \:direction   #1\p@
                      \:NormalizeDeg\:direction  }
\Define\:RotateTo(1){ \:direction  #1\p@
                      \:NormalizeDeg\:direction  }

\Define\:tRotate(2){
   \advance\:Vdirection  #2\p@
   \:NormalizeDeg\:Vdirection  \:Rotate(#1)}
\Define\:tRotateTo(2){
   \:Vdirection #2\p@
   \:NormalizeDeg\:Vdirection  \:RotateTo(#1)}
\Define\LineF(1){\MarkLoc(,)\MoveF(#1) {\LineToLoc(,)}}
\Define\:MvF(1){    \:CosSin\:direction
   \:d#1\:Xunits   \advance\:X  \:Cons\:d\:x
   \:d#1\:Yunits   \advance\:Y  \:Cons\:d\:y   }
\Define\MoveFToOval(2){
   \:NormalizeDeg\:direction \:CosSin\:direction  \:Zunits\p@
   \:d#2\:Zunits    \:x\:Cons\:d\:x
   \ifdim \:d=\z@  \:err5{...,\:Cons\:d} \fi
   \:d#1\:Zunits    \:y\:Cons\:d\:y
   \ifdim \:d=\z@  \:err5{\:Cons\:d,...} \fi     \:SearchDir
   \XSaveUnits
      \if:rotated
        \:Zunits\p@     \:Zunits#1\:Zunits
        \:x\:Cons\:Zunits\:Xunits
        \:Zunits\p@     \:Zunits#2\:Zunits
        \:y\:Cons\:Zunits\:Yunits
      \else
        \:x#1\:Xunits  \:y#2\:Yunits  \fi     \relax
      \Units(\:x,\:y)
      \:Vdirection\:direction  \:direction\:ddd
      \MoveF(\@ne)  \:direction\:Vdirection  \XRecallUnits  }
\Define\:tMvF(1){    \:CosSin\:Vdirection
   \:d #1\:Zunits   \advance\:Z  \:Cons\:y\:d
   \:xx#1\:Xunits     \:yy#1\:Yunits
   \:xx\:Cons\:x\:xx    \:yy\:Cons\:x\:yy
   \:CosSin\:direction
   \:x\:Cons\:xx\:x    \:y\:Cons\:yy\:y
   \advance\:X  \:x   \advance\:Y  \:y } \Define\CSeg{\:Opt[]\:CSeg1}
\def\:CSeg[#1]#2(#3,#4){   \MarkLoc($^)
   \MoveToLoc(#4) \:x\:X \:y\:Y
   \if:IIID  \:d\:Z  \fi   \MoveToLoc(#3)
   \advance\:x -\:X  \:x#1\:x
   \advance\:y -\:Y  \:y#1\:y
   \if:IIID   \advance\:d -\:Z  \:d#1\:d \fi
   \:t{#2}
   \edef\:temp{\the\:t(
      \expandafter\:Frac\the\:x\noexpand\:du,
      \expandafter\:Frac\the\:y\noexpand\:du \if:IIID ,
      \expandafter\:Frac\the\:d\noexpand\:du \fi)}
   \MoveToLoc($^)     \:temp}\Define\LSeg{\:Opt[]\:LSeg1}
\def\:LSeg[#1]#2(#3,#4){   \:distance(#3,#4)
   \:d#1\:d  \:t{#2}
   \edef\:temp{\the\:t(\expandafter\:Frac\the\:d\noexpand\:du)}
   \:temp}\Define\DSeg{\:Opt[]\:DSeg1}

\def\:DSeg[#1]#2(#3,#4){   \MarkLoc(^)
   \MoveToLoc(#4)  \:xxx\:X   \:yyy\:Y \:xxxx\:Z
   \MoveToLoc(#3)
   \advance\:xxx -\:X   \advance\:yyy -\:Y
   \ifdim \:AbsVal\:xxx<\:mmmp
      \ifdim \:AbsVal\:yyy<\:mmmp \:wrn5{#3,#4}
   \fi\fi
   \if:IIID
      \advance \:xxxx  -\:Z
      \:divide\:xxx\:Xunitsx
      \:divide\:yyy\:Yunitsy
      \:divide\:xxxx\:Zunits
      \:x\:xxx  \:y\:yyy   \:ddd\z@  \:length
      \:x\:d    \:y\:xxxx  \:SearchDir
      \:yyyy\:ddd
      \:x\:xxx  \:y\:yyy
   \else
      \:x  \:AbsVal\:Yunitsy
\:y  \:AbsVal\:Xunitsy  \ifdim \:y>\:x  \:x\:y  \fi
\:y  \:AbsVal\:Yunitsx  \ifdim \:y>\:x  \:x\:y  \fi
\:y  \:AbsVal\:Xunitsx  \ifdim \:y>\:x  \:x\:y  \fi
\:K  \:InCons\:x  \relax
      \ifnum \:K<\thr@@     \:K\@ne
\else \ifnum \:K<\sixt@@n   \:K4
\else \ifnum \:K<\:XC       \:K\sixt@@n
\else \ifnum \:K<\@m        \:K\@cclvi
\fi \fi \fi \fi
\divide\:xxx \:K
\divide\:yyy \:K  
      \:x \:Cons\:Yunitsy\:xxx
\advance\:x -\:Cons\:Yunitsx\:yyy
      \:y-\:Cons\:Xunitsy\:xxx
\advance\:y \:Cons\:Xunitsx\:yyy     
   \fi
   \:SearchDir   \:ddd#1\:ddd   \:t{#2}
   \edef\:temp{\the\:t(
      \:Cons\:ddd \if:IIID ,\:Cons\:yyyy \fi)}
   \MoveToLoc(^) \:temp}   \def\:theDoReg{\def\DoReg{\the\:DoReg}}

\Define\MarkLoc{  \:theDoReg
   \expandafter\edef \csname \:MarkLoc}

\Define\MarkGLoc{  \:theDoReg
   \expandafter\xdef \csname \:MarkLoc}

\Define\:MarkLoc(1){ Loc\space#1:\endcsname{
   \:X\the\:X  \:Y\the\:Y
   \if:IIID \:Z\the\:Z \fi  }       \let\DoReg\:DoReg }

\Define\MoveToLoc(1){   \:theDoReg   \expandafter\ifx
   \csname Loc\space#1:\endcsname\relax \:err2{#1}\fi
   \csname Loc\space#1:\endcsname     \let\DoReg\:DoReg  }

\Define\MarkPLoc(1){  \:theDoReg
   \if:IIID   \:project
      \expandafter\edef \csname Loc\space#1:\endcsname{
         \:X\the\:x  \:Y\the\:y  \:Z\z@}
   \else \:err1\MarkPLoc  \fi   \let\DoReg\:DoReg}

\Define\WriteLoc(1){{ \:theDoReg \edef\:temp{#1}
   \ifx \:temp\empty \else \MoveToLoc(#1) \fi
   \immediate\write\sixt@@n{...
      \:temp=(\the\:X,\the\:Y\if:IIID,\the\:Z\fi)}}}

\Define\:LnToLoc(1){
   \:x\:X \:y\:Y   \MoveToLoc(#1)
   { \:LnTo(\:x\du,\:y\du) }}

\Define\:tLnToLoc(1){
   \:xx\:X \:yy\:Y \:dd\:Z     \MoveToLoc(#1)
   { \:tLnTo(\:xx\du,\:yy\du,\:dd\du) }}\def\:GetLine#1#2#3#4#5{
   \MoveToLoc(#1)   \divide\:X \:eight  \divide\:Y \:eight
   #3\:X  #4\:Y
   \MoveToLoc(#2)   \divide\:X \:eight  \divide\:Y \:eight
   \advance #3 -\:X  \advance #4 -\:Y
   #5\:Cons#3\:Y      \advance #5 -\:Cons#4\:X
   \divide #3 \:eight     \divide  #4 \:eight \relax      }
\def\MoveToLL(#1,#2)(#3,#4){
   \:GetLine{#1}{#2}\:x \:y \:xxx
   \:GetLine{#3}{#4}\:xx\:yy\:xxxx
   \:ddd \:Cons\:x \:yy     \advance\:ddd -\:Cons\:xx\:y
   \ifdim  \:AbsVal\:ddd < \:mmmp
      \:X\@cclv\p@  \:Y\:X
      \:wrn3{(\string#1,\string#2)(\string#3,\string#4)}
   \else
      \:divide\:xxx\:ddd    \:divide\:xxxx\:ddd
      \:X\:Cons\:xxx\:xx   \advance\:X -\:Cons\:xxxx\:x
      \:Y\:Cons\:xxx\:yy   \advance\:Y -\:Cons\:xxxx\:y
   \fi   }\Define\MoveToCC{\:Opt[]\:MoveToCC{}}

\def\:MoveToCC[#1](#2,#3)(#4,#5){
  \:UserUnits(#2,#3)(#4,#5)
  \:distance($#2,$#4)
\ifnum \:d<\:mp \MoveToLoc(#3)
   \:wrn3{(#1,#2)(#3,#4)}
\else                 \:xx \:d
   \:distance($#2,$#3)  \:xxx\:d
   \:distance($#5,$#4) 
     \:yy \:xxx  \advance\:yy -\:d
\:yyy\:xxx  \advance\:yyy \:d
\:divide\:yy\:xx     \:yy\:Cons\:yy\:yyy
\advance\:yy \:xx  \divide\:yy \tw@
     \:yyy \ifdim \:AbsVal\:xxx>\:AbsVal\:yy \:xxx \else \:yy \fi
\ifdim \:AbsVal\:yyy<\:mp \:yyy\z@ \else
   \:divide\:xxx\:yyy  \:sqr\:xxx
   \:yyyy\:yy  \:divide\:yyyy\:yyy  \:sqr\:yyyy
   \advance\:xxx -\:yyyy   \:Sqrt\:xxx
   \:yyy\:Cons\:yyy\:xxx
\fi 
     \MoveToLoc($#4)  \:x\:X  \:y\:Y  \:divide\:yy\:xx
\MoveToLoc($#2)
\advance\:x -\:X     \advance\:y -\:Y
\advance\:X \:Cons\:yy\:x
\advance\:Y \:Cons\:yy\:y
     \:divide\:yyy\:xx
\advance\:X  #1\:Cons\:yyy\:y
\advance\:Y -#1\:Cons\:yyy\:x

  \fi  \:SysUnits  }\def\:UserUnits(#1,#2)(#3,#4){
   \:xx\:Xunitsx  \:xxx\:Xunitsy
   \:yy\:Yunitsx  \:yyy\:Yunitsy
   \:xxxx\:Cons\:yy\:xxx  \advance\:xxxx \:Cons\:yyy\:xx
   \ifdim \:AbsVal\:xxxx>\:mmmp
      \:divide\:xx\:xxxx   \:divide\:xxx{-\:xxxx}
      \:divide\:yy{\:xxxx} \:divide\:yyy\:xxxx
   \fi
   \:UnLoc(#1)  \:UnLoc(#2)
   \:UnLoc(#3)  \:UnLoc(#4)}\Define\:SysUnits{
   \MoveTo(\:Cons\:X,\:Cons\:Y)}\Define\:UnLoc(1){
   \MoveToLoc(#1)  \:d\:Cons\:yyy\:X
   \advance\:d \:Cons\:yy\:Y
   \:Y\:Cons\:xx\:Y
   \advance\:Y  \:Cons\:xxx\:X   \:X\:d
   \MarkLoc($#1)}

\def\:MoveToLC[#1](#2,#3)(#4,#5){
   \:UserUnits(#2,#3)(#4,#5)
   \MoveToLoc($#2)  \:x\:X  \:y\:Y
\MoveToLoc($#3)  \advance\:x -\:X
                \advance\:y -\:Y
   \edef\:temp{ \:xxx\the\:x  \:yyy\the\:y }
\MoveToLoc($#4)
\advance\:X \:y  \advance\:Y -\:x
\MarkLoc(^$) \MoveToLL($#4,^$)($#2,$#3)
\MarkLoc(^$)  
   \:distance($#4,$#5)  \:xx\:d  \:distance($#4,^$)
\ifdim      \:d>\:xx        \:wrn3{(#2,#3)(#4,#5)}
\else \ifdim \:d<\:mmp     \:yy\:xx   \else   \:yy\:d
      \:divide\:yy\:xx  \:sqr\:yy
      \:yy-\:yy   \advance\:yy \p@
      \:Sqrt\:yy  \:yy\:Cons\:xx\:yy
\fi \fi
   \:temp   \:x\:xxx  \:y\:yyy  \:length  \:xx\:d
\:divide\:xxx\:xx
\:divide\:yyy\:xx
\advance\:X  #1\:Cons\:yy\:xxx
\advance\:Y  #1\:Cons\:yy\:yyy
   \:SysUnits }  \def\Object#1{\:Opt(){\:DefineSD#1}0}

\def\:DefineSD#1(#2){\begingroup  \:multid#1
   \:DraCatCodes   \:DefSD#1(#2)}

\def\:DefSD#1(#2)#3{
   \expandafter\::Define\csname\string#1.\endcsname(#2){
      \:t{\:SubD{#3}}
      \if:IIID \edef\:temp{\noexpand\TwoDim \the\:t
                           \noexpand\EndTwoDim}
      \else    \def\:temp{\the\:t}   \fi         \:temp}
   \def#1{\def\:SDname{\csname\string#1.\endcsname}
          \:Opt[]\:CallSD{}}}

\def\:CallSD[#1]{ \edef\:Entry{#1} \:SDname }\def\:SubD#1{
   \let\::RecallXLoc\:AddXLoc   \gdef\:AddXLoc{}
   \edef\:RecallBor{ \global\:LBorder  \the\:LBorder
                  \global\:RBorder  \the\:RBorder
                  \global\:TeXLoc\the\:TeXLoc }
\global\:TeXLoc\z@
\setbox\:box\vbox{\EntryExit(0,0,0,0)
\begingroup
   
   \:InDraw  #1
\endgroup
\:SetDrawWidth                    \let\:XLoc\relax
\xdef\:AddXLoc{\:dd\the\:LBorder  \:AddXLoc}}
\:RecallBor
   \:ddd\dp\:box
   \ifx \:Entry\empty
   \:DrawBox
\else
   \let\:RecallIn\:InOut
   \:x\:X    \:y\:Y
   \def\:XLoc(##1,##2,##3){
      \def\:temp{##1}
      \ifx \:temp\:Entry \:X\:x  \advance\:X -##2
                         \:Y\:y  \advance\:Y -##3
      \fi}
   \:AddXLoc
   \advance\:X  \:dd      \advance\:Y -\:ddd
   \EntryExit(-1,-1,\:InOut2,\:InOut3)  \:DrawBox
   \let\:InOut\:RecallIn
\fi
   \MarkLoc(^)
      \MoveToExit(-1,-1)
   \:xxx\:X   \:yyy\:Y   \advance\:yyy  \:ddd
   \def\:XLoc(##1,##2,##3){
      \:X\:xxx  \advance\:X  ##2     \advance\:X -\:dd
      \:Y\:yyy  \advance\:Y  ##3
      \MarkLoc(##1)}
   \:AddXLoc
   \MoveToLoc(^)
   \ifx \:Entry\empty \else     \MoveToLoc(\:Entry) \fi
   \global\let\:AddXLoc\::RecallXLoc }
\Define\:MarkXLoc(1){  \:theDoReg
   \let\:XLoc\relax
   \xdef\:AddXLoc{\:AddXLoc \:XLoc(#1,\the\:X,\the\:Y)}
   \let\DoReg\:DoReg} 
   \catcode`\ 10 \catcode`\^^M5 \catcode`\^^I10
   \let\wlog\:wlog  \let\:wlog\:undefined
\:RestoreCatcodes     \tracingstats1   

\newcommand{\VolumeHeader}{}
\newcommand{\VolumeSerial}{LNS}

\newcommand{\ActivityName}{ {\normalsize {\it 
Summer School on High-dimensional Manifold Topology }}}
\newcommand{\ActivityDate}{ {\normalsize {\it
Trieste, 21 May -- 8 June 2001
}}}

\usepackage{amssymb}
\newcommand{\be}{\begin{equation}}
\newcommand{\ee}{\end{equation}}
\newcommand{\bea}{\begin{eqnarray}}
\newcommand{\eea}{\end{eqnarray}}

\def\K{\mathbb K}

\def\Q{\mathbb Q}
\def\R{\mathbb R}
\def\Z{\mathbb Z}
\def\GT{\mathcal G \mathcal T}
\def\GCCT{{\mathcal G \mathfrak C \mathcal C \mathcal T}}
\def\qed{$\square$}
\newcommand{\LectureHeader}{Circle valued Morse theory}

\begin{document}
\pagestyle{myheadings}
\markboth{\LectureHeader}{\VolumeHeader}
\markright{\VolumeHeader}


\begin{titlepage}


\title{Circle valued Morse theory and Novikov homology}

\author{Andrew Ranicki\thanks{aar@maths.ed.ac.uk}
\\[1cm]
{\normalsize
{\it Department of Mathematics and Statistics}}\\
{\normalsize
{\it University of Edinburgh, Scotland, UK} }
\\[10cm]
{\normalsize {\it Lecture given at the: }}
\\
\ActivityName 
\\
\ActivityDate 
\\[1cm]
{\small \VolumeSerial} 
}
\date{}
\maketitle
\thispagestyle{empty}
\end{titlepage}

\baselineskip=14pt
\newpage
\thispagestyle{empty}

\begin{abstract}

Traditional Morse theory deals with real valued functions $f:M\to \R$
and ordinary homology $H_*(M)$.  The critical points of a Morse
function $f$ generate the Morse-Smale complex $C^{MS}(f)$ over $\Z$,
using the gradient flow to define the differentials.  The isomorphism
$H_*(C^{MS}(f))\cong H_*(M)$ imposes homological restrictions on real
valued Morse functions.  There is also a universal coefficient version
of the Morse-Smale complex, involving the universal cover
$\widetilde{M}$ and the fundamental group ring $\Z[\pi_1(M)]$.

The more recent Morse theory of circle valued functions $f:M \to S^1$
is more complicated, but shares many features of the real valued
theory.  The critical points of a Morse function $f$ generate the
Novikov complex $C^{Nov}(f)$ over the Novikov ring $\Z((z))$ of formal
power series with integer coefficients, using the gradient flow of the
real valued Morse function $\overline{f}:\overline{M}=f^*\R \to \R$ on
the infinite cyclic cover to define the differentials.  The Novikov
homology $H^{Nov}_*(M)$ is the $\Z((z))$-coefficient homology of
$\overline{M}$.  The isomorphism $H_*(C^{Nov}(f))\cong H^{Nov}_*(M)$
imposes homological restrictions on circle valued Morse functions.

Chapter 1 reviews real valued Morse theory.  Chapters 2,3,4 introduce
circle valued Morse theory and the universal coefficient versions of
the Novikov complex and Novikov homology, which involve the universal
cover $\widetilde{M}$ and a completion $\widehat{\Z[\pi_1(M)]}$ of
$\Z[\pi_1(M)]$.  Chapter 5 formulates an algebraic chain complex model 
(in the universal coefficient version) for 
the relationship between the $\Z((z))$-module Novikov complex $C^{Nov}(f)$ 
of a circle valued Morse function $f:M \to S^1$ and the $\Z$-module Morse-Smale
complex $C^{MS}(f_N)$ of the real valued Morse function
$f_N=\overline{f}\vert:M_N=\overline{f}^{-1}[0,1] \to [0,1]$ on a
fundamental domain of the infinite cyclic cover $\overline{M}$.

\end{abstract}
\vfill


{\it AMS numbers:} 57R70, 55U15

\newpage
\thispagestyle{empty}
\tableofcontents

\newpage
\setcounter{page}{1}

\section{Introduction}

The Morse theory of circle valued functions $f:M \to S^1$
relates the topology of a manifold $M$ to the critical points of $f$,
generalizing the traditional theory of real valued Morse functions 
$M \to \R$. However, the relationship is somewhat more complicated
in the circle valued case than in the real valued case, and the
roles of the fundamental group $\pi_1(M)$ and of the choice of 
gradient-like vector field $v$ are more significant (and less
well understood).

The Morse-Smale complex $C=C^{MS}(M,f,v)$ is defined geometrically for
a real valued Morse function $f:M^m \to \R$ and a suitable choice of
gradient-like vector field $v:M \to \tau_M$.  In general, there is a
$\Z[\pi]$-coefficient Morse-Smale complex for each group morphism
$\pi_1(M) \to \pi$, with
$$C_i~=~\Z[\pi]^{c_i(f)}$$
if there are $c_i(f)$ critical points of index $i$.
The differentials $d:C_i \to C_{i-1}$ are
defined by counting the $\widetilde{v}$-gradient flow lines in the
cover $\widetilde{M}$ of $M$ classified by $\pi_1(M) \to \pi$.
In the simplest case $\pi=\{1\}$ this is just $\widetilde{M}=M$, 
and if $p \in M$ is a
critical point of index $i$ and $q \in M$ is a critical point of index
$i-1$ the $(p,q)$-coefficient in $d$ is the number $n(p,q)$ of lines
from $p$ to $q$, with sign chosen according to orientations.
The homology of the Morse-Smale complex is 
isomorphic to the ordinary homology of $M$
$$H_*(C^{MS}(M,f,v))~\cong~H_*(M)$$ 
so that
\begin{itemize}
\item[(a)] the critical points of $f$ can be used to compute $H_*(M)$,
\item[(b)] $H_*(M)$ provides lower bounds on the number 
of critical points in any Morse function $f:M \to \R$, which must have 
at least as many critical points of index $i$ as there are $\Z$-module
generators for $H_i(M)$ (Morse inequalities).
\end{itemize}
Basic real valued Morse theory is reviewed in Chapter 2.

In the last 40 years there has been much interest in the Morse theory
of circle valued functions $f:M^m \to S^1$, starting with the work of
Stallings \cite{Stallings}, Browder and Levine \cite{BrowderLevine},
Farrell \cite{Farrell} and Siebenmann \cite{Siebenmann} on the
characterization of the maps $f$ which are homotopic to the
projections of fibre bundles over $S^1$ : these are the circle valued
Morse functions without any critical points.

About 20 years ago, Novikov 
(\cite{Novikov1981},\cite{Novikov1982},\cite{Novikov1993},\cite{Novikov1996} 
(pp. 194--199)) was motivated by problems in physics and dynamical
systems to initiate the general Morse theory of closed 1-forms,
including circle valued functions $f:M \to S^1$ as the most important
special case. The new idea was to use the {\it Novikov ring} of formal
power series with an infinite number of positive coefficients and
a finite number of negative coefficients
$$\Z((z))~=~\Z[[z]][z^{-1}]~=\{
\sum\limits^{\infty}_{j=-\infty}n_jz^j\,\vert\,
n_j \in \Z,~\hbox{$n_j=0$ for all $j<k$, for some $k$}\}$$
as a counting device for the gradient flow lines of the
real valued Morse function $\overline{f}:\overline{M}=f^*\R \to \R$ on
the (non-compact) infinite cyclic cover $\overline{M}$ of $M$,
with the indeterminate $z$ corresponding to the generating covering
translation $z:\overline{M} \to\overline{M}$. 
For $f$ the number of gradient flow lines starting at a critical point
$p \in M$ is finite in the generic case. On the other hand,
for $\overline{f}$ the 
number of gradient flow lines starting at a critical point
$\overline{p} \in \overline{M}$ may be infinite in the generic case,
so the counting methods for real and circle valued Morse theory 
are necessarily different.

The {\it Novikov complex} $\widehat{C}=C^{Nov}(M,f,v)$ is defined for a
circle valued Morse function $f:M^m \to S^1$ and a suitable choice of
gradient-like vector field $v:M \to \tau_M$.  In general, there is a
$\widehat{\Z[\Pi]}$-coefficient Novikov complex for each factorization
of $f_*:\pi_1(M) \to \pi_1(S^1)=\Z$ as $\pi_1(M) \to \Pi \to \Z$, with
$\widehat{\Z[\Pi]}$ a completion of $\Z[\Pi]$, with
$$\widehat{C}_i~=~\widehat{\Z[\Pi]}^{c_i(f)}$$
if there are $c_i(f)$ critical points of index $i$.
The differentials $d:C_i \to C_{i-1}$ are
defined by counting the $\widetilde{v}$-gradient flow lines in the
cover $\widetilde{M}$ of $M$ classified by $\pi_1(M) \to \Pi$.
The construction of the Novikov complex for arbitrary
$\widehat{Z[\Pi]}$ is described in Chapter 3. In the simplest case
$$\Pi~=~\Z~,~\Z[\Pi]~=~\Z[z,z^{-1}]~,~\widehat{\Z[\Pi]}~=~\Z((z))~,~
\widetilde{M}~=~\overline{M}~=~f^*\R~.$$ 
For a critical point $\overline{p} \in \overline{M}$ of index $i$
and a critical point $\overline{q} \in \overline{M}$ of an index $i-1$ 
the $(\overline{p},\overline{q})$-coefficients in $\widehat{d}$ is
$$\widehat{n}(\overline{p},\overline{q})~=~\sum\limits^{\infty}_{j=k}
n(\overline{p},z^j\overline{q})z^j\in \Z((z))$$
with $n(\overline{p},z^j\overline{q})$ the signed number of
$\overline{v}$-gradient flow lines of the real valued Morse function
$\overline{f}:\overline{M} \to \R$ from $\overline{p}$ to the translate
$z^j\overline{q}$ of $\overline{q}$, and 
$k=[\overline{f}(\overline{p}) - \overline{f}(\overline{q})]$.
The convention is that the generating covering translation
$z:\overline{M} \to \overline{M}$ is to be chosen 
parallel to the downward gradient flow $v:M \to \tau_M$, with 
$$\overline{f}(zx)~=~\overline{f}(x)-1 \in \R ~~(x \in \overline{M})~.$$
In particular, this means that for $f=1:M=S^1\to S^1$
$$z~:~\overline{M}~=~\R \to \overline{M}~=~\R~;~x \mapsto x-1~.$$
\indent
Circle valued Morse theory is necessarily more complicated than real
valued Morse theory.  The Morse-Smale complex $C^{MS}(M,f:M \to \R,v)$
is an absolute object, describing $M$ on the chain level,
with $c_0(f)>0$, $c_m(f)>0$.  This is the
algebraic analogue of the fact that every continuous function $f:M \to
\R$ on a compact space attains an absolute minimum and an absolute
maximum. By contrast, the Novikov complex $C^{Nov}(M,f:M \to S^1,v)$
is a relative object, measuring the chain level difference between $f$
and the projection of a fibre bundle (= Morse function with no critical
points).  A continuous function $f:M \to S^1$ can just go round and
round! The connection between the geometric properties of $f$ and the
algebraic topology of $M$ is still not yet completely understood,
although there has been much progress in the work of Pajitnov, Farber,
the author and others.

The {\it Novikov homology groups} of a space $M$ with respect to a
cohomology class $f \in [M,S^1]=H^1(M)$ are defined by
$$H^{Nov}_*(M,f)~=~H_*\big(\Z((z))\otimes_{\Z[z,z^{-1}]}C(\overline{M})\big)~.$$
The homology groups of the Novikov complex are isomorphic to the
Novikov homology groups
$$H_*(C^{Nov}(M,f,v))~\cong~H^{Nov}_*(M,f)~.$$
By analogy with the real valued case :
\begin{itemize}
\item[(a)] the critical points of $f$ can be used to compute
$H^{Nov}_*(M,f)$,
\item[(b)] $H^{Nov}_*(M,f)$ provides lower bounds on the number of critical
points in any Morse function $f:M \to S^1$, which must have at least as
many critical points of index $i$ as there are 
$\Z((z))$-module generators for $H^{Nov}_i(M,f)$ (Morse-Novikov 
inequalities).
\end{itemize}
Novikov homology is constructed in Chapter 4, for arbitrary
$\widehat{\Z[\Pi]}$-coefficients..

Novikov conjectured (\cite{Arnold}) that for a generic class of
gradient-like vector fields $v \in \GT(f)$ the functions $j \mapsto
n(\overline{p},z^j\overline{q})$ have subexponential growth.  

Let $S \subset \Z[z]$ be the subring of the polynomials $s(z)$ such that
$s(0)=1$, which (up to sign) are precisely the polynomials 
which are invertible in the power series ring $\Z[[z]]$.
The localization $S^{-1}\Z[z,z^{-1}]$ of $\Z[z,z^{-1}]$ is identified
with the subring of  $\Z((z))$ consisting of the quotients
$\displaystyle{r(z) \over s(z)}$ with $r(z) \in \Z[z,z^{-1}]$, $s(z) \in S$.

Pajitnov \cite{Pajitnov1997},\cite{Pajitnov1999} constructed a $C^0$-dense
subspace $\GCCT(f) \subset \GT(f)$ of gradient-like vector fields $v$
for which the differentials in the Novikov complex $C^{Nov}(M,f,v)$ are
rational
$$\widehat{n}(\overline{p},\overline{q})~=~
\sum\limits^{\infty}_{j=-\infty}n(\overline{p},z^j\overline{q})z^j \in
S^{-1}\Z[z,z^{-1}] \subset \Z((z))$$
and the functions
$j \mapsto n(\overline{p},z^j\overline{q})$ have polynomial growth.
The idea is to cut $M$ along the inverse image $N=f^{-1}(0)$
(assuming $0 \in S^1$ is a regular value of $f$), giving a 
fundamental domain 
$$(M_N;N,z^{-1}N)~=~\overline{f}^{-1}([0,1];\{0\},\{1\})$$
for $\overline{f}:\overline{M} \to \R$, and to then use a kind of
cellular approximation theorem to give a chain level approximation to
the gradient flow in
$$(f_N,v_N)~=~(\overline{f},\overline{v})\vert~:
(M_N,f_N,v_N)\to ([0,1];\{0\},\{1\})~.$$
The mechanism described in Chapter 5 below then gives a chain complex 
over $S^{-1}\Z[z,z^{-1}]$ inducing $C^{Nov}(M,f,v)$. 
Hutchings and Lee \cite{HutchingsLee1},\cite{HutchingsLee2} used a
similar method to get enough information from $C^{Nov}(M,f,v)$ for
generic $v$ to obtain an estimate on the number of closed $v$-gradient
flow lines $\gamma:S^1 \to M$. 

Farber and Ranicki \cite{FarberRanicki} and Ranicki \cite{Ranicki1999}
constructed an `algebraic Novikov complex' in $S^{-1}\Z[z,z^{-1}]$ for
any circle Morse valued function $f:M \to S^1$, using any $CW$
structure on $N=f^{-1}(0)$, the extension to a $CW$ structure on $M_N$,
and a cellular approximation to the inclusion $z^{-1}N \to M_N$.  The
construction is recalled in Chapter 5, including the non simply
connected version.  In many cases (e.g.  for $v \in \GCCT(f)$) this
algebraic model does actually coincide with the geometric Novikov
complex $C^{Nov}(M,f,v)$.

The Morse-Novikov theory of circle valued functions on
finite-dimensional manifolds and Novikov homology have many
applications to symplectic topology, Floer homology, and Seiberg-Witten
theory (Po\'zniak \cite{Pozniak}, Le and Ono \cite{LeOno}, Hutchings
and Lee \cite{HutchingsLee1}, \cite{HutchingsLee2}, $\dots$).  Also,
circle valued Morse theory on infinite-dimensional manifolds features
in the work of Taubes on Casson's homology 3-sphere invariant and gauge
theory.  However, these notes are not a survey of all the applications
of circle valued Morse theory and Novikov homology! They deal
exclusively with the basic development in the finite-dimensional case
and some of the applications to the classification of manifolds.

\section{Real valued Morse theory}

This section reviews the real valued Morse theory, which is a
prerequisite for circle valued Morse theory.  The traditional
references Milnor \cite{Milnor1963}, \cite{Milnor1965} remain the best
introductions to real valued Morse theory.  Bott \cite{Bott} gives a
beautiful account of the history of Morse theory, including the
development of the modern chain complex point of view inspired by
Witten.

Let $M$ be a compact differentiable $m$-dimensional manifold.
The {\it critical points} of a differentiable function $f:M \to \R$ are
the zeros $p \in M$ of the differential $\nabla f : \tau_M \to \tau_{\R}$.
A {\it Morse function} $f:M \to \R$ is a differentiable function in which 
every critical point $p \in M$ is required to be isolated and nondegenerate,
meaning that in local coordinates
$$f(p+(x_1,x_2,\dots,x_m))~=~
f(p)-\sum\limits^i_{j=1}(x_j)^2+\sum\limits^m_{j=i+1}(x_j)^2$$
with $i$ the index of $p$. 
The subspace of Morse functions is $C^2$-dense in the space of all 
differentiable functions $f:M \to \R$.

A vector field $v:M \to \tau_M$ is {\it gradient-like} for $f$ if
there exists a Riemannian metric $\langle~,~\rangle$
on $M$ such that  
$$\langle v,w \rangle~=~\nabla f(w) \in \R~~(w \in \tau_M)~.$$
Note that $\langle~,~\rangle$ and $\nabla f$ determine $v$, and that
the zeros of $v$ are the critical points of $f$.

A {\it $v$-gradient flow line} $\gamma: \R \to M$ satisfies
$$\gamma'(t)~=~-v(\gamma(t))  \in \tau_M(\gamma(t))~~(t \in \R)~.$$ 
The minus sign here gives the {\it downward} gradient flow, so that 
$$f(\gamma(s)) > f(\gamma(t))~{\rm if}~s<t~.$$
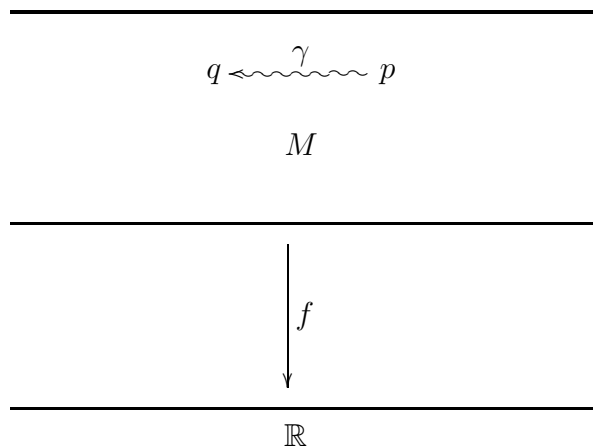
$$\Draw
\LineAt(-20,0,200,0)
\LineAt(-20,80,200,80)
\MoveTo(90,60)
\Text(--$\xymatrix@R-12pt { q && \ar@{~>}[ll]_-{\displaystyle{\gamma}} p}$--)
\MoveTo(90,30)
\Text(--$M$--)
\MoveTo(100,-35)
\Text(--$\xymatrix@R+30pt{\ar[d]^-{\displaystyle{f}} & \\ & }$--)
\LineAt(-20,-70,200,-70)
\MoveTo(88,-80)
\Text(--$\R$--)
\EndDraw$$
\smallskip

\noindent The limits
$$\lim_{t \to -\infty} \gamma(t)~=~p~~,~~\lim_{t \to \infty}
\gamma(t)~=~q \in M$$
are critical points of $f$ with $f(q)<f(p)$, and if $\gamma$ is isolated
then
$${\rm index}(q)~=~{\rm index}(p)-1~.$$
For every point $x \in M$ there is a $v$-gradient flow line
$\gamma_x:\R \to M$ (which is unique up to scaling) such that
$\gamma_x(0)=x \in M$. If $x$ is a critical point take $\gamma_x$ to
be the constant path at $x$.

The {\it unstable} and {\it stable} manifolds of a critical point $p
\in M$ of index $i$ are the open manifolds
$$W^u(p,v)~=~\{x\in M\,\vert\,\lim_{t\rightarrow -\infty}\gamma_x(t)=p\}~,~
W^s(p,v) ~=~\{x\in M\,\vert\,\lim_{t\rightarrow\infty}\gamma_x(t)=p\}~.$$
The unstable and stable manifolds are images of immersions 
$\R^i \to M$, $\R^{m-i}\to M$ respectively, which are embeddings
near $p \in M$.

The basic results relating a Morse function $f:M^m \to \R$ to the 
topology of $M$ concern the inverse images
$$N_a~=~f^{-1}(a)$$
of the regular values $a \in \R$, which are closed $(m-1)$-dimensional
manifolds, and the cobordisms
$$(M_{a,b};N_a,N_b)~=~f^{-1}([a,b];\{a\},\{b\})~~(a <b)~.$$

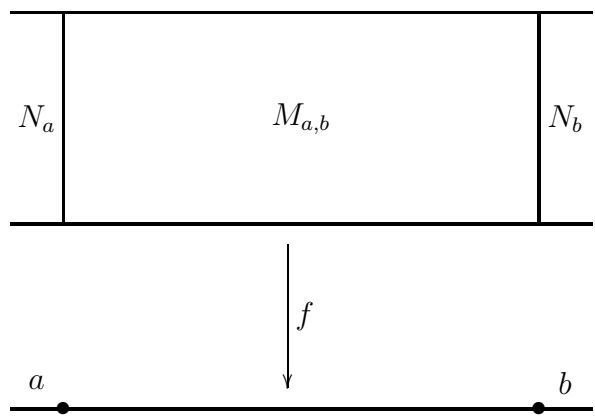
$$\Draw
\DrawRectAt(0,0,180,80)
\LineAt(-20,0,0,0)
\LineAt(180,0,200,0)
\LineAt(-20,80,0,80)
\LineAt(180,80,200,80)
\MoveTo(-10,40)
\Text(--$N_a$--)
\MoveTo(190,40)
\Text(--$N_b$--)
\MoveTo(90,40)
\Text(--$M_{a,b}$--)
\MoveTo(100,-35)
\Text(--$\xymatrix@R+30pt{\ar[d]^-{\displaystyle{f}} & \\ & }$--)
\LineAt(-20,-70,200,-70)
\MoveTo(0,-70)
\Text(--$\bullet$--)
\MoveTo(180,-70)
\Text(--$\bullet$--)
\MoveTo(-10,-60)
\Text(--$a$--)
\MoveTo(190,-60)
\Text(--$b$--)
\EndDraw$$

\noindent The results are:
\begin{itemize}
\item[(i)] if $[a,b] \subset \R$ contains no critical values
the $v$-gradient flow determines a diffeomorphism 
$$N_b \to N_a~;~x \mapsto \gamma_x((f\gamma_x)^{-1}(a))~,$$
\item[(ii)] if $[a,b] \subset \R$ contains a unique critical
value $f(p) \in (a,b)$, and $p \in M$ is a critical point of index $i$,
then $N_b$ is obtained from $N_a$ by surgery on a tubular neighbourhood
$S^{i-1} \times D^{m-i} \subset N_a$ of $S^{i-1}=W^u(p,v) \cap N_a$
$$N_b~=~N_a \backslash (S^{i-1} \times D^{m-i}) \cup D^i \times S^{m-i-1}$$
with $D^i \times S^{m-i-1} \subset N_b$ a tubular neighbourhood of
$S^{m-i-1}=W^s(p,v) \cap N_b$, and $(M_{a,b};N_a,N_b)$ the trace of the surgery
$$M_{a,b}~=~N_a \times [0,1] \cup D^i \times D^{m-i}~.$$
\end{itemize}

Let $\GT(f)$ denote the set of gradient-like vector fields $v$ on $M$
which satisfy the Morse-Smale transversality condition that for any
critical points $p,q \in M$ with ${\rm index}(p)=i$, ${\rm index}(q)=j$ 
the submanifolds $W^u(p,v)^i$, $W^s(q,v)^{m-j} \subset M^m$
intersect transversely in an $(i-j)$-dimensional submanifold
$W^u(p,v)\cap W^s(q,v) \subset M$.
The subspace $\GT(f)$ is dense in the space of gradient-like vector fields for $f$.

Suppose that the Morse function $f:M \to \R$ has $c_i(f)$ critical
points of $f$ of index $i$, and that the critical points
$p_0,p_1,p_2,\dots \in M$ are arranged to satisfy
$${\rm index}(p_0) \leqslant {\rm index}(p_1) \leqslant {\rm index}(p_2) \leqslant \dots~,~
f(p_0) <f(p_1) < f(p_2) < \dots~.$$
A choice of $v \in \GT(f)$ determines a handle decomposition of $M$
$$M~=~\bigcup\limits_{i=0}^m\bigcup_{c_i(f)} D^i \times D^{m-i}$$
with one $i$-handle $h^i=D^i \times D^{m-i}$ for each critical point of
index $i$.

The {\it Morse-Smale complex} $C^{MS}(M,f,v)$ is defined for
a Morse-Smale pair
$(f:M \to \R,v\in \GT(f))$ and a regular cover
$\widetilde{M}$ of $M$ with group of covering translations $\pi$,
to be the based f.g.  free $\Z[\pi]$-module chain complex with
$$d~:~C^{MS}(M,f,v)_i~=~\Z[\pi]^{c_i(f)} \to 
C^{MS}(M,f,v)_{i-1}~=~\Z[\pi]^{c_{i-1}(f)}~;~ 
\widetilde{p} \mapsto \sum\limits_{\widetilde{q}}
n(\widetilde{p},\widetilde{q})\widetilde{q}$$
with $n(\widetilde{p},\widetilde{q})\in \Z$ the finite signed
number of $\widetilde{v}$-gradient flow lines 
$\widetilde{\gamma}:\R \to \widetilde{M}$
which start at a critical point $\widetilde{p} \in
\widetilde{M}$ of $\widetilde{f}:\widetilde{M}\to \R$ with index $i$ and
terminate at a critical point $\widetilde{q} \in \widetilde{M}$ of
index $i-1$.  Choose an arbitrary lift of each critical point
$p \in M$ of $f$ to a critical point $\widetilde{p} \in \widetilde{M}$
of $\widetilde{f}$, obtaining a basis for $C^{MS}(M,f,v)$.
The Morse-Smale complex is the cellular chain complex
$$C^{MS}(M,f,v)~=~C(\widetilde{M})$$
of the $CW$ structure on $\widetilde{M}$ in which the $i$-cells are the
lifts of the $i$-handles $h^i$. In particular, the homology of the
Morse-Smale complex is the ordinary homology of $\widetilde{M}$
$$H_*(C^{MS}(M,f,v))~=~H_*(\widetilde{M})~.$$
\indent If $(f,v):M \to \R$ is modified to $(f',v'):M \to \R$ by adding
a pair of critical points $p,q$ of index $i,i-1$ with $n(p,q,v)=1$ the
Morse-Smale complex $C^{MS}(M,f',v')$ is obtained from $C^{MS}(M,f,v)$
by attaching an elementary chain complex
$$E~:~\dots \to 0 \to E_i~=~\Z[\pi] \xrightarrow[]{1} E_{i-1}~=~\Z[\pi] \to 
0 \to \dots~,$$
with an exact sequence
$$0 \to C^{MS}(M,f,v) \to C^{MS}(M,f',v') \to E \to 0 ~.$$
Conversely, if $m \geqslant 5$ then the Whitney trick applies to realize
the elementary moves of Whitehead torsion theory by cancellation of
pairs of critical points (or equivalently, handles).
This cancellation is the basis of the proofs of the $h$- and $s$-cobordism
theorems.

The identity $C^{MS}(M,f,v)=C(M)$ (for $\widetilde{M}=M$)
gives the {\it Morse inequalities}
$$c_i(f) \geqslant b_i(M) + q_i(M)+q_{i-1}(M)$$
with
$$b_i(M)~=~{\rm dim}_{\Z}\big(H_i(M)/T_i(M))~,~
q_i(M)~=~\#\,T_i(M)$$
the Betti numbers of $M$, where 
$$T_i(M)~=~\{x \in H_i(M)\,\vert\,nx=0~\text{for some}~n \neq 0 \in \Z\}$$
is the torsion subgroup of $H_i(M)$ and $\#$ denotes the minimum number
of generators.
Smale used the cancellation of critical points to prove that these
inequalities are sharp for $\pi_1(M)=\{1\}$, $m \geqslant 5$: there exists
$(f,v):M \to \R$ with the minimum possible number of critical points
$$c_i(f)~=~b_i(M)+q_i(M)+q_{i-1}(M)~.$$
The method is to start with an arbitrary Morse function $f:M \to \R$,
and to systematically cancel pairs of critical points until this is no 
longer possible. 

The Morse-Smale complex $C^{MS}(M,f,v)$ is also defined for a
Morse function on an $m$-dimensional cobordism
$f:(M;N,N') \to ([0,1];\{0\},\{1\})$
with $v \in \GT(f)$. In this case there is a relative handle decomposition
$$M~=~N \times [0,1] \cup \bigcup\limits^m_{i=0}\bigcup\limits_{c_i(f)}D^i
\times D^{m-i}$$
and $C^{MS}(M,f,v)=C(\widetilde{M},\widetilde{N})$. The $s$-cobordism
theorem states that for a Morse function $f$ on an $h$-cobordism 
$\tau(C^{MS}(M,f,v))=0 \in Wh(\pi_1(M))$ if (and for $m \geqslant 6$ only if)
the critical points of $f$ can be stably cancelled in pairs.

\section{The Novikov complex}

Morse functions $f:M \to S^1$, gradient-like vector field $v$, critical
points, index, $c_i(f)$, are defined in the same way as for the real
valued case in Chapter 1.  Again, the subspace of Morse functions is
$C^2$-dense in the space of all functions $f:M \to S^1$. But it is
harder to decide which pairs of critical points can be cancelled.

A Morse function $f:M \to S^1$ lifts to a $\Z$-equivariant Morse
function $\overline{f}:\overline{M}=f^*\R \to \R$ on the infinite
cyclic cover 
$$\xymatrix@R+18pt@C+18pt{
\overline{M} \ar[r] \ar[d]_-{\displaystyle{\overline{f}}} 
& M \ar[d]^-{\displaystyle{f}} \\
\R \ar[r] & S^1}$$
\smallskip

\noindent Let $z:\overline{M} \to \overline{M}$ be the downward
generating covering translation. 
$$\Draw
\LineAt(-20,0,200,0)
\LineAt(-20,80,200,80)
\MoveTo(90,60)
\Text(--$\xymatrix@R-12pt { && \ar@{~>}[ll]_-{\displaystyle{z}}}$--)
\MoveTo(90,30)
\Text(--$\overline{M}$--)
\MoveTo(103,-35)
\Text(--$\xymatrix@R+30pt{\ar[d]^-{\displaystyle{\overline{f}}} & \\ & }$--)
\LineAt(-20,-70,200,-70)
\MoveTo(88,-80)
\Text(--$\R$--)
\EndDraw$$
\smallskip

\noindent Let $\GT(f)$ be the space of gradient-like vector fields $v:M
\to \tau_M$ such that a lift $\overline{v}:\overline{M} \to
\tau_{\overline{M}}$ satisfies the Morse-Smale transversality
condition.  The Novikov complex of a circle valued Morse function is
defined by analogy with the Morse-Smale complex of a real valued
function, as follows.

Given a ring $A$ and an automorphism $\alpha:A \to A$ let $z$
be an indeterminate over $A$ with
$$az~=~z\alpha(a)~~~(a \in A)~.$$
\indent The {\it $\alpha$-twisted Laurent polynomial extension} of $A$ is
the localization of the $\alpha$-twisted polynomial extension 
$A_{\alpha}[z]$ inverting $z$
$$A_{\alpha}[z,z^{-1}]~=~A_{\alpha}[z][z^{-1}]~,$$
the ring of  polynomials 
$\sum\limits_{j=-\infty}^{\infty} a_jz^j$ ($a_j \in A$) 
such that $\{j \in \Z \,\vert\,a_j \neq 0\}$ is finite.

The {\it $\alpha$-twisted Novikov ring} of $A$ is the
localization of the completion of $A_{\alpha}[z]$
$$A_{\alpha}((z))~=~A_{\alpha}[[z]][z^{-1}]~,$$
the ring of power series 
$\sum\limits_{j=-\infty}^{\infty} a_jz^j$ ($a_j \in A$) 
such that $\{j \leqslant 0 \,\vert\,a_j \neq 0\}$ is finite.

Given $f:M \to S^1$ let $\widetilde{M}$ be a regular cover of
$\overline{M}$, with group of covering translations $\pi$.  Only the
case of connected $M,\overline{M},\widetilde{M}$ will be considered. 
Let $\Pi$ be the group of covering translations of $\widetilde{M}$ over
$M$, so that there is defined a group extension
$$\{1\} \to \pi \to \Pi \to \Z \to \{1\}$$
with a lift of $1 \in \Z$ to an element $z \in \Pi$
such that the covering translation $z:\widetilde{M} \to \widetilde{M}$
induces $z:\overline{M} \to \overline{M}$ on $\overline{M}=\widetilde{M}/\pi$.
Thus 
$$\Pi~=~\pi\times_{\alpha}\Z~~,~~\Z[\Pi]~=~\Z[\pi]_{\alpha}[z,z^{-1}]~.$$
Write the $\alpha$-twisted Novikov ring as
$$\widehat{\Z[\Pi]}~=~\Z[\pi]_{\alpha}((z))~.$$
Choose a lift of each critical point $p \in M$ of $f$ to a critical
point $\widetilde{p} \in \widetilde{M}$ of $\widetilde{f}$.

The {\it Novikov complex} $C^{Nov}(M,f,v)$ of 
$(f:M \to S^1,v\in \GT(f))$ is the based
f.g.  free $\widehat{\Z[\Pi]}$-module chain complex with
$$\begin{array}{l}
d~:~C^{Nov}(M,f,v)_i~=~\Z[\pi]_{\alpha}((z))^{c_i(f)} \to 
C^{Nov}(M,f,v)_{i-1}~=~\Z[\pi]_{\alpha}((z))^{c_{i-1}(f)}~;\\
\hskip150pt \widetilde{p} \mapsto \sum\limits^
{\infty}_{j=-\infty}\sum\limits_{\widetilde{q}}
n(\widetilde{p},z^j\widetilde{q}) z^j\widetilde{q}
\end{array}$$
with $n(\widetilde{p},\widetilde{q})\in \Z$ the finite signed number of 
$\widetilde{v}$-gradient flow lines $\widetilde{\gamma}:\R \to \widetilde{M}$ 
which start at a critical point $\widetilde{p} \in
\widetilde{M}$ of $\widetilde{f}:\widetilde{M}\to \R$ with index $i$ and
terminate at a critical point $\widetilde{q} \in \widetilde{M}$ of index $i-1$.  

\noindent{\it Exercise.} Work out $C^{Nov}(S^1,f,v)$ for
$$f~:~S^1 \to S^1~;~[t] \mapsto [4t-9t^2+6t^3]~~(0 \leqslant t \leqslant 1)~.\eqno{\square}$$

The original definition of Novikov
\cite{Novikov1981},\cite{Novikov1982} was in the special case
$$\widetilde{M}~=~\overline{M}~,~\pi~=~\{1\}~,~\Pi~=~\Z~,~\alpha~=~1~,~
\widehat{\Z[\Pi]}~=~\Z((z))$$ 
when $C^{Nov}(M,f,v)$ is a based f.g. 
free $\Z((z))$-module chain complex (as in the Exercise).

Take $\widetilde{M}$ to be the universal cover of $M$ and
$\pi=\pi_1(\overline{M})$, $\alpha:\pi \to \pi$ the automorphism
induced by a generating covering translation $z:\overline{M} \to
\overline{M}$, $\Pi=\pi_1(M)=\pi\times_{\alpha}\Z$.  This case gives
the based f.g.  free $\widehat{\Z[\pi_1(M)]}$-module Novikov complex
$C^{Nov}(M,f,v)$ of Pajitnov \cite{Pajitnov1996}.

There is only one class of Morse functions $f:M \to S^1$ for which the
Novikov complex is easy to compute:

\noindent{\it Example.} Let $M$ be the mapping torus of a diffeomorphism
$h:N\to N$ of a closed $(m-1)$-dimensional manifold
$$M~=~T(h)~=~(N \times [0,1])/\{(x,0) \sim (h(x),1)\}~.$$
$$\Draw
\Text(--$T(h)$--)
\DrawOvalArc(25,25)(-45,225)
\DrawOvalArc(75,75)(-45,225)
\MoveTo(35.36,-35.36)
\Text(--$\bullet$--)
\MoveTo(-35.36,-35.36)
\Text(--$\bullet$--)
\MoveTo(0,-50)
\Text(--$(x,0) \sim (h(x),1)$--)
\LineAt(17.68,-17.68,53.04,-53.04)
\LineAt(-17.68,-17.68,-53.04,-53.04)
\MoveTo(0,50)
\Text(--$N\times [0,1]$--)
\EndDraw$$
\medskip

\noindent The fibre bundle projection
$$f~:~M~=~T(h) \to S^1~=~[0,1]/\{0 \sim 1\}~;~[x,t] \mapsto [t]$$
has no critical points, so that $C^{Nov}(M,f,v)=0$ for any $v \in \GT(f)$.
\hfill\qed

\section{Novikov homology}

The Novikov homology $H_*^{Nov}(M,f;\widehat{\Z[\Pi]})$ is defined
for a space $M$ with a map $f:M \to S^1$ and a factorization
of $f_*:\pi_1(M) \to \pi_1(S^1)$ through a group $\Pi$. 
The relevance of the Novikov complex $C^{Nov}(M,f,v)$ to the Morse
theory of a Morse map $f:M \to S^1$ is immediately obvious.  
The relevance of the Novikov homology is rather less obvious,
even though there are isomorphisms 
$H_*(C^{Nov}(M,f,v))\cong H_*^{Nov}(M,f;\widehat{\Z[\Pi]})$ !

The $R$-coefficient homology of a space $M$ is
defined for any ring morphism $\Z[\pi_1(M)] \to R$ 
$$H_*(M;R)~=~H_*(C(M;R))$$
using any free $\Z[\pi_1(M)]$-module chain complex $C(\widetilde{M})$ 
(e.g. cellular, if $M$ is a $CW$ complex) and
$C(M;R)=R\otimes_{\Z[\pi_1(M)]}C(\widetilde{M})$.

Given a group $\pi$ and an automorphism $\alpha:\pi \to \pi$ let
$\pi\times_{\alpha}\Z$ be the group with elements $gz^j$ ($g \in \pi$, 
$j \in \Z$), and multiplication by $gz=\alpha(g)z$, so that
$$\Z[\pi\times_{\alpha}\Z]~=~\Z[\pi]_{\alpha}[z,z^{-1}]~.$$
For any map $f:M \to S^1$ with $M$ connected the infinite cyclic
cover $\overline{M}=f^*\R$ is connected if and only if 
$f_*:\pi_1(M) \to \pi_1(S^1)=\Z$ is onto, in which case
$$\pi_1(M)~=~\pi_1(\overline{M})\times_{\alpha_M}\Z$$
with $\alpha_M:\pi_1(\overline{M}) \to \pi_1(\overline{M})$ the
automorphism induced by a generating covering translation 
$z:\overline{M} \to \overline{M}$.

Suppose given a connected space $M$ with a cohomology class 
$f \in [M,S^1]=H^1(M)$ such that $\overline{M}=f^*\R$ is connected.
Given a factorization of the surjection $f_*:\pi_1(M) \to \pi_1(S^1)$ 
$$f_*~:~\pi_1(M)~=~\pi_1(\overline{M})\times_{\alpha_M}\Z \to \Pi \to \Z$$
let $\pi=\hbox{\rm ker}(\Pi \to \Z)$ and let $z \in \Pi$ be the image
of $z=(0,1) \in \pi_1(M)$, so that $\Pi=\pi\times_{\alpha}\Z$ with
$$\alpha~:~\pi \to \pi~;~g \mapsto z^{-1}gz~.$$ 
The
$\widehat{\Z[\Pi]}$-coefficient {\it Novikov homology} of $(M,f)$ is 
$$H^{Nov}_*(M,f;\widehat{\Z[\Pi]})~=~H_*(M;\widehat{\Z[\Pi]})~,$$
with $\widehat{\Z[\Pi]}=\Z[\pi]_{\alpha}((z))$.

In the original case 
$$\widetilde{M}~=~\overline{M}~,~\pi~=~\{1\}~,~\Pi~=~\Z~,~
\widehat{\Z[\Pi]}~=~\Z((z))~,$$
and $H^{Nov}_*(M,f;\widehat{\Z[\Pi]})$
may be written as $H^{Nov}_*(M,f)$, or even just $H^{Nov}_*(M)$.

\noindent{\it Example 1.} 
The $\Z((z))$-coefficient cellular chain complex of $S^1$ is
$$C\big(S^1;\Z((z))\big)~:~\dots \to 0 \to \Z((z)) \xrightarrow[]{1-z} \Z((z))$$
and $1-z \in \Z((z))$ is a unit, so  $H^{Nov}_*(S^1)=0$.
\hfill\qed

\noindent{\it Example 2.} Let $N$ be a connected finite $CW$ complex with
cellular $\Z$-module chain complex $C(N)$, and let $h:N \to N$ be a self-map 
with induced $\Z$-module chain map $h:C(N) \to C(N)$.
The $\Z((z))$-coefficient cellular chain complex of the mapping torus $T(h)$
with respect to the canonical projection  
$$f~:~T(h)\to S^1~;~[x,t] \mapsto [t]$$
is the algebraic mapping cone
$$C^+\big(T(h);\Z((z))\big)~=~{\mathcal C}\big(1-zh:C(N)((z)) \to C(N)((z))\big)~.$$
Now $1-zh$ is a $\Z((z))$-module chain equivalence, so that
$$H^{Nov}_*(T(h),f)~=~0~.$$ 
The $\Z((z))$-coefficient cellular chain complex of the mapping torus $T(h)$
with respect to the other projection  
$$-f~:~T(h)\to S^1~;~[x,t] \mapsto [1-t]$$
is the algebraic mapping cone
$$C^-\big(T(h);\Z((z))\big)~=~{\mathcal C}\big(z-h:C(N)((z)) \to C(N)((z))\big)~.$$
If $h:N \to N$ is a homotopy equivalence then $z-h$  is a
$\Z((z))$-module chain equivalence, so that
$$H^{Nov}_*(T(h),-f)~=~0~,$$
but in general $H^{Nov}_*(T(h),-f)\neq 0$ -- see Example 3 below
for an explicit non-zero example.\\
\hbox to\hsize{\hfill\qed}

\noindent{\it Example 3.} 
The Novikov homology groups of the mapping torus $T(2)$ of 
the double covering map $2:S^1 \to S^1$ are
$$\begin{array}{l} 
H^{Nov}_1(T(2),f)~=~\Z((z))/(1-2z)~=~0~,\\[1ex]
H^{Nov}_1(T(2),-f)~=~\Z((z))/(z-2)~=~\widehat{\Z}_2[1/2]~=~\widehat{\Q}_2
\neq 0
\end{array}$$
with $\widehat{\Q}_2$ the 2-adic field
(Example 23.25 of Hughes and Ranicki \cite{HughesRanicki}). 
The inverse of 
$$n~=~2^a(2b+1)\in \Z$$
is
$$n^{-1}~=~z^{-a}(1-zb+z^2b^2-z^3b^3+\dots) \in
\Z((z))/(2-z)~=~\widehat{\Q}_2~.\eqno{\hbox{\qed}}$$

\noindent{\it Theorem.} (Novikov \cite{Novikov1981}, \cite{Novikov1982}
for $\pi=\{1\}$, Pajitnov \cite{Pajitnov1995}) \\
{\it The Novikov complex $C^{Nov}(M,f,v)$ is $\widehat{\Z[\Pi]}$-module 
chain equivalent to $C(M;\widehat{\Z[\Pi]})$, with isomorphisms}
$$H_*(C^{Nov}(M,f,v))~\cong~H^{Nov}_*(M,f;\widehat{\Z[\Pi}])~.
\eqno{\hbox{\qed}}$$

The chain equivalence $C^{Nov}(M,f,v)\simeq C(M;\widehat{\Z[\Pi}])$
will be described in Chapter 4 below.

The Novikov ring $\Z((z))$ is a principal ideal domain, and
$H^{Nov}_*(M,f)$ is the homology of a f.g.  free $\Z((z))$-module chain
complex.  Thus each $H^{Nov}_i(M,f)$ is a f.g.  $\Z((z))$-module, which
splits as free$\oplus$torsion, by the structure theorem for f.g. 
modules over a principal ideal domain.

The {\it Novikov numbers} of any finite $CW$ complex $M$ and
$f \in H^1(M)$ are the Betti numbers of Novikov homology
$$b^{Nov}_i(M,f)~=~{\rm dim}_{\Z((z))}\big(H^{Nov}_i(M,f)/T^{Nov}_i(M,f)\big)~,~
q^{Nov}_i(M,f)~=~\#\,T^{Nov}_i(M,f)$$
where  
$$T^{Nov}_i(M,f)~=~
\{x \in H^{Nov}_i(M,f)\,\vert\,ax=0~{\rm for~some}~a \neq 0 \in \Z((z))\}$$
is the torsion $\Z((z))$-submodule of $H^{Nov}_i(M,f)$,
and $\#$ denotes the minimum number of generators.

The {\it Morse-Novikov inequalities} (\cite{Novikov1981})
$$c_i(f) \geqslant b^{Nov}_i(M,f) + q^{Nov}_i(M,f)+q^{Nov}_{i-1}(M,f)$$
are an immediate consequence of the isomorphisms
$H_*(C^{Nov}(M,f,v))\cong H^{Nov}_*(M,f)$, since for any 
f.g. free chain complex $C$ over a principal ideal domain $R$
$${\rm dim}_R(C_i)\geqslant b_i(C)+q_i(C)+q_{i-1}(C)$$
where
$$b_i(C)~=~{\rm dim}_R\big(H_i(C)/T_i(C)\big)~,~
q_i(C)~=~\#\,T_i(C)$$
with 
$$T_i(C)~=~\{x \in H_i(C)\,\vert\,rx=0~{\rm for~some}~r \neq 0  \in R\}$$
the $R$-torsion submodule of $H_i(C)$, and $\#$ denoting the minimal
number of $R$-module generators.

Farber \cite{Farber1985} proved that the Morse-Novikov inequalities are
sharp for $\pi_1(M)=\Z$, $m \geqslant 6$ : for any such manifold there
exists a Morse function $f:M \to S^1$ representing $1 \in [M,S^1]=H^1(M)$
with the minimum possible numbers of critical points
$$c_i(f)~=~b^{Nov}_i(M,f) + q^{Nov}_i(M,f)+q^{Nov}_{i-1}(M,f)~.$$
Again, the method is to start with an arbitrary Morse function $f:M \to S^1$
in the homotopy class, and to systematically cancel pairs of
critical points until this is no longer possible. 

When does the Novikov homology vanish?

\noindent{\it Proposition} (Ranicki \cite{Ranicki1995})
Let $A$ be a ring with an automorphism $\alpha:A \to A$.
A finite f.g. free $A_{\alpha}[z,z^{-1}]$-module chain complex $C$
is such that
$$H_*(A_{\alpha}((z))\otimes_{A_{\alpha}[z,z^{-1}]}C)~=~
H_*(A_{\alpha}((z^{-1}))\otimes_{A_{\alpha}[z,z^{-1}]}C)~=~0$$
if and only if $C$ is $A$-module chain equivalent to a finite
f.g. projective $A$-module chain complex.\hfill\qed

Note that for an algebraic Poincar\'e complex $(C,\phi)$ 
$$H_*(A_{\alpha}((z))\otimes_{A_{\alpha}[z,z^{-1}]}C)~=~0~
\hbox{if and only if}~
H_*(A_{\alpha}((z^{-1}))\otimes_{A_{\alpha}[z,z^{-1}]}C)~=~0~,$$
so the two Novikov homology vanishing conditions can be replaced by just one.

Recall that a space $X$ is {\it finitely dominated} if there exist a
finite $CW$ complex and maps $i:X \to K$, $j:K \to X$ such that $ji
\simeq 1:X \to X$.  Wall \cite{Wall} proved that a $CW$ complex
$X$ is finitely dominated if and only if $\pi_1(X)$ is finitely
presented and the cellular chain complex $C(\widetilde{X})$ of the
universal cover $\widetilde{X}$ is chain equivalent to a finite f.g. 
projective $\Z[\pi_1(X)]$-module chain complex.

In the simply-connected case $\pi_1(\overline{M})=\{1\}$ the 
following conditions on a map $f:M \to S^1$ from an $m$-dimensional manifold
$M$ are equivalent :
\begin{itemize}
\item[(i)] $\overline{M}$ is finitely dominated,
\item[(ii)] $\overline{M}$ is homotopy equivalent to a finite $CW$ complex,
\item[(iii)] $H^{Nov}_*(M,f)=0$,
\item[(iv)] $b^{Nov}_i(M,f)=q^{Nov}_i(M,f)=0$,
\item[(v)] $C(\overline{M})$ is chain equivalent to a finite f.g.
free $\Z$-module chain complex,
\item[(vi)] the homology groups $H_*(\overline{M})$ are f.g. $\Z$-modules.
\end{itemize}
Browder and Levine \cite{BrowderLevine} used handle exchanges (= the
ambient surgery version of the cancellation of adjacent critical
points) to prove that (vi) holds if (and for $m \geqslant 6$ only if) $f:M
\to S^1$ is homotopic to the projection of a fibre bundle.

Farrell \cite{Farrell} and Siebenmann \cite{Siebenmann} defined
a Whitehead torsion obstruction $\Phi(M,f) \in Wh(\pi_1(M))$
for a map $f:M^m \to S^1$ with finitely dominated $\overline{M}=f^*\R$,
such that $\Phi(M,f)=0$ if (and for $m \geqslant 6$ only if) $f$ is
homotopic to the projection of a fibre bundle.

\noindent{\it Theorem} (Ranicki \cite{Ranicki1995})\\
(i) {\it For any finite $CW$ complex $M$ and map $f:M \to S^1$ 
the infinite cyclic cover $\overline{M}=f^*\R$ of $M$ is finitely dominated if and only if
$\pi_1(\overline{M})$ is finitely presented and
$$H^{Nov}_*(M,f;\widehat{\Z[\pi_1(M)]})~=~0~.$$}
(ii) {\it For any Morse map $f:M \to S^1$ on an $m$-dimensional manifold $M$
with finitely dominated $\overline{M}$ the torsion of the Novikov complex 
$\tau(C^{Nov}(M,f,v)) \in K_1(\widehat{\Z[\pi_1(M)]})/I$
determines and is determined by the Farrell-Siebenmann fibering obstruction 
$\Phi(M,f) \in Wh(\pi_1(M))$, where $I \subseteq K_1(\widehat{\Z[\pi_1(M)]})$ is the subgroup
generated by $\pm \pi_1(M)$ and $\tau(1-zh)$ for square matrices $h$
over $\Z[\pi_1(\overline{M})]$. Thus $\tau(C^{Nov}(M,f,v)) \in I$
if {\rm (}and for $m \geqslant 6$ only if\,{\rm )}
$f$ is homotopic to a fibre bundle.}\hfill\qed

See Chapter 22 of Hughes and Ranicki \cite{HughesRanicki} and Chapter 15
of Ranicki \cite{Ranicki1998} for more detailed accounts of the
relationship between the torsion of the Novikov complex and the
Farrell-Siebenmann fibering obstruction.

See Latour \cite{Latour} and Pajitnov \cite{Pajitnov1996} for 
circle Morse-theoretic
proofs that if $m\geqslant 6$, $\overline{M}$ is finitely dominated
and $\tau(C^{Nov}(M,f,v))\in I$ then it is possible to pairwise 
cancel all the critical points of $f$.

\section{The algebraic model for circle valued Morse theory}

In many cases the Novikov complex $C^{Nov}(M,f,\alpha)$ of a
circle valued Morse function $f:M \to S^1$ can be constructed from an
algebraic model for the $\overline{v}$-gradient flow in a fundamental
domain of the infinite cyclic cover $\overline{M}$. 

An {\it algebraic fundamental domain} $(D,E,F,g,h)$
consists of finite based f.g. free $A$-module chain complexes
$D,E$ and chain maps $g:D \to E$, $h:z^{-1}D \to E$ of the form
$$\begin{array}{l}
d_E~=~\begin{pmatrix} d_D & c \\ 0 & d_F \end{pmatrix}~:~
E_i~=~D_i \oplus F_i \to E_{i-1}~=~D_{i-1} \oplus F_{i-1}~,\\[2ex]
g~=~\begin{pmatrix} 1 \\ 0 \end{pmatrix}~:~
D_i \to E_i~=~D_i \oplus F_i~,\\[2ex]
h~=~\begin{pmatrix} h_D \\ h_F \end{pmatrix}~:~
z^{-1}D_i \to E_i~=~D_i \oplus F_i~.
\end{array}$$
$$\Draw
\LineAt(-20,0,320,0)
\LineAt(-20,80,320,80)
\LineAt(0,0,0,80)
\LineAt(100,0,100,80)
\LineAt(200,0,200,80)
\LineAt(300,0,300,80)
\MoveTo(-10,40)
\Text(--$zD$--)
\MoveTo(50,40)
\Text(--$zE$--)
\MoveTo(90,40)
\Text(--$D$--)
\MoveTo(122,45)
\Text(--$\xymatrix{\ar[r]^-{\displaystyle{g}}&}$--)
\MoveTo(150,40)
\Text(--$E$--)
\MoveTo(178,45)
\Text(--$\xymatrix{&\ar[l]_-{\displaystyle{zh}}}$--)
\MoveTo(220,40)
\Text(--$z^{-1}D$--)
\MoveTo(260,40)
\Text(--$z^{-1}E$--)
\MoveTo(320,40)
\Text(--$z^{-2}D$--)
\EndDraw$$

\noindent
Define the {\it algebraic Novikov complex} $\widehat{F}$ to be the
based f.g. free $A_{\alpha}((z))$-module chain complex with
$$\begin{array}{ll}
d_{\widehat{F}}&=~d_F+zh_F(1-zh_D)^{-1}c\\[1ex]
&=~d_F+\sum\limits^{\infty}_{j=1}z^jh_F(h_D)^{j-1}c~:~
\widehat{F}_i~=~(F_i)_{\alpha}((z)) \to
\widehat{F}_{i-1}~=~(F_{i-1})_{\alpha}((z))~,
\end{array}$$
as in Farber and Ranicki \cite{FarberRanicki} and Ranicki \cite{Ranicki1999}.
The $A_{\alpha}((z))$-module chain map
$$\phi~=~g-zh~=~\begin{pmatrix} 1-zh_D \\-zh_F \end{pmatrix}~
:~D_{\alpha}((z)) \to E_{\alpha}((z))$$
is a split injection in each degree (since $1-zh_D$ is an isomorphism), 
and the inclusions $F_i \to E_i$ determine a canonical isomorphism
of based f.g. free $A_{\alpha}[z,z^{-1}]$-module chain complexes
$$\widehat{F}~\cong~{\rm coker}(\phi)~.$$
\indent Here is how algebraic fundamental domains and the algebraic
Novikov complex arise in topology.

Let $f:M \to S^1$ be a Morse function with regular value $0 \in S^1$. 
$$\Draw
\DrawCircle(25)
\DrawCircle(75)
\LineAt(0,-25,0,-75)
\MoveTo(10,-50)
\Text(--$N$--)
\MoveTo(0,50)
\Text(--$M_N$--)
\MoveTo(122,0)
\Text(--$
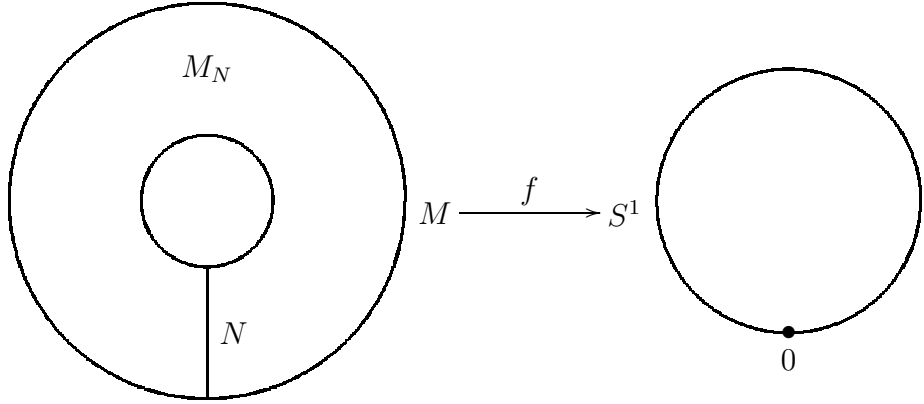
\xymatrix@C+29pt{M\ar[r]^-{\displaystyle{f}}&S^1}$--)
\MoveTo(220,0)
\DrawCircle(50)
\MoveTo(220,-50)
\Text(--$\bullet$--)
\MoveTo(220,-60)
\Text(--$0$--)
\EndDraw$$
\smallskip

\noindent Cut $M$ along the inverse image
$$N^{m-1}~=~f^{-1}(0) \subset M$$ 
to obtain a geometric fundamental domain 
$$(M_N;N,z^{-1}N)~=~\overline{f}^{-1}([0,1];\{0\},\{1\})$$
for the infinite cyclic cover 
$$\overline{M}~=~f^*\R~=~\bigcup\limits^{\infty}_{j=-\infty}z^jM_N~.$$
The restriction 
$$f_N~=~\overline{f}\vert~:~(M_N;N,z^{-1}N) \to ([0,1];\{0\},\{1\})$$
is a real valued Morse function with $v_N=\overline{v}\vert \in \GT(f_N)$.
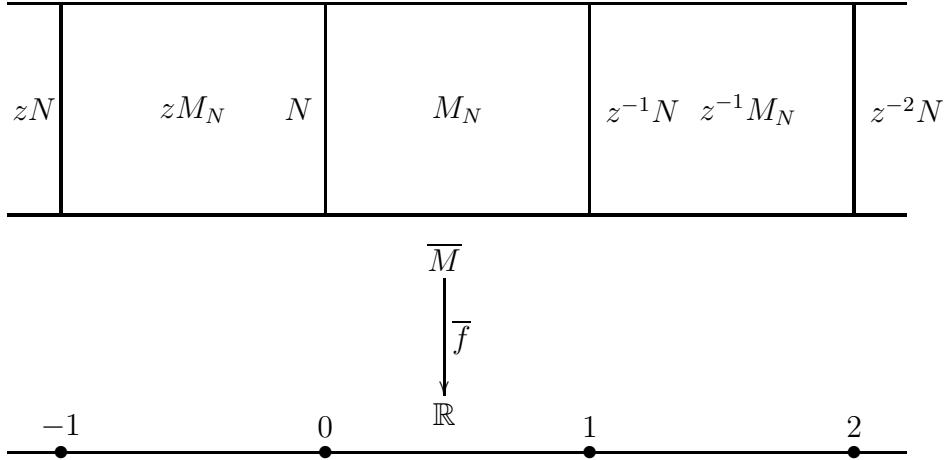
$$\Draw
\LineAt(-20,0,320,0)
\LineAt(-20,80,320,80)
\LineAt(0,0,0,80)
\LineAt(100,0,100,80)
\LineAt(200,0,200,80)
\LineAt(300,0,300,80)
\MoveTo(-10,40)
\Text(--$zN$--)
\MoveTo(50,40)
\Text(--$zM_N$--)
\MoveTo(90,40)
\Text(--$N$--)
\MoveTo(150,40)
\Text(--$M_N$--)
\MoveTo(220,40)
\Text(--$z^{-1}N$--)
\MoveTo(260,40)
\Text(--$z^{-1}M_N$--)
\MoveTo(320,40)
\Text(--$z^{-2}N$--)
\MoveTo(160,-45)
\Text(--$\xymatrix@R+20pt{\overline{M}\ar[d]^-{\displaystyle{\overline{f}}} & \\ \R& }$--)
\LineAt(-20,-90,320,-90)
\MoveTo(0,-90)
\Text(--$\bullet$--)
\MoveTo(100,-90)
\Text(--$\bullet$--)
\MoveTo(200,-90)
\Text(--$\bullet$--)
\MoveTo(300,-90)
\Text(--$\bullet$--)
\MoveTo(0,-80)
\Text(--$-1$--)
\MoveTo(100,-80)
\Text(--$0$--)
\MoveTo(200,-80)
\Text(--$1$--)
\MoveTo(300,-80)
\Text(--$2$--)
\EndDraw$$
\smallskip

\noindent The cobordism $(M_N;N,z^{-1}N)$ has a handlebody decomposition
$$M_N~=~N \times [0,1] \cup
\bigcup\limits^m_{i=0}\bigcup\limits_{c_i(f)}D^i \times D^{m-i}$$
with one $i$-handle for each index $i$ critical point of $f$.  Given a
$CW$ structure on $N$ with $c_i(N)$ $i$-cells use this handlebody
decomposition to define a $CW$ structure on $M_N$ with
$c_i(N)+c_i(f)$ $i$-cells. A regular cover $\widetilde{M}$
of $\overline{M}$ with group of covering translations $\pi$
is a regular cover of $M$ with group of covering translations 
$\Pi=\pi\times_{\alpha}\Z$ (as before), with
$$\Z[\Pi]~=~\Z[\pi]_{\alpha}[z,z^{-1}]~~,~~\widehat{\Z[\Pi]}~=~
\Z[\Pi]_{\alpha}((z))~.$$
Use a cellular approximation $h:z^{-1}N \to M_N$ to the inclusion to
define an algebraic fundamental domain $(D,E,F,g,h)$ over $A=\Z[\pi]$
$$D~=~C(\widetilde{N})~~,~~E~=~C(\widetilde{M}_N)~~,~~
F~=~C^{MS}(M_N,f_N,v_N)~=~C(\widetilde{M}_N,\widetilde{N})~.$$
The mapping cylinder of $h:N \to M_N$ is a $CW$ complex $M'_N$ with
two copies of $N$ as subcomplexes. Identifying these copies there
is obtained a $CW$ complex structure on $M$ with
$\widehat{\Z[\Pi]}$-coefficient cellular chain complex
$$C(M;\widehat{\Z[\Pi]})~=~{\mathcal C}(\phi)$$
the algebraic mapping cone of the $\widehat{\Z[\Pi]}$-module chain map
$$\phi~=~g-zh~:~D_{\alpha}((z)) \to E_{\alpha}((z))~,$$
with
$$\begin{array}{l}
d_{{\mathcal C}(\phi)}~=~\begin{pmatrix}
-d_D & 0 & 0 \\
1-zh_D & d_D & c \\
-zh_F & 0 & d_D \end{pmatrix}~:\\[4ex]
{\mathcal C}(\phi)_i~=~(D_{i-1}\oplus D_i \oplus F_i)_{\alpha}[z,z^{-1}] \to 
{\mathcal C}(\phi)_{i-1}~=~(D_{i-2}\oplus D_{i-1} \oplus F_{i-1})_{\alpha}[z,z^{-1}]~.
\end{array}$$
\noindent The algebraic Novikov complex $\widehat{F}={\rm coker}(\phi)$ 
is a based f.g. free $\widehat{\Z[\Pi]}$-module chain complex such that
$${\rm dim}_{\widehat{\Z[\Pi]}}\widehat{F}_i~=~c_i(f)~.$$
In many cases $\widehat{F}=C^{Nov}(M,f,v)$, and in even more cases
$\widehat{F}$ is simple isomorphic to $C^{Nov}(M,f,v)$.
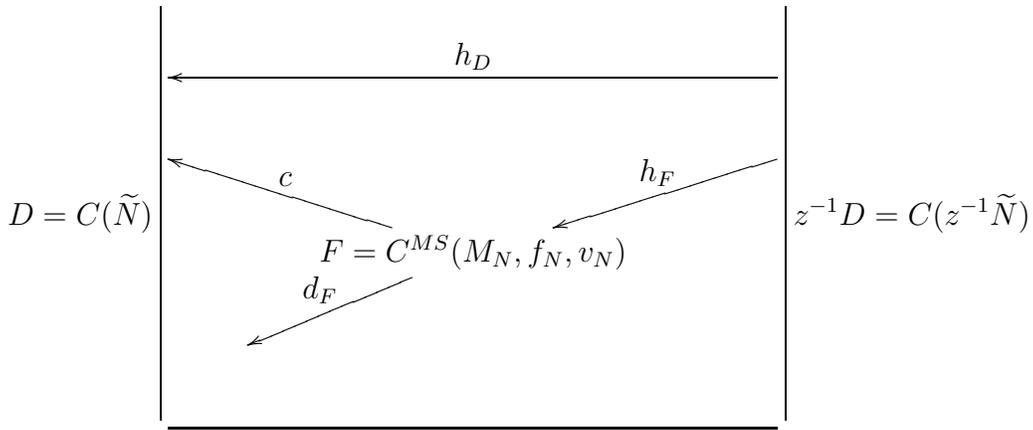
$$\xymatrix
{ \ar@{-}[rrrr]
\ar@{-}[ddddd]_-{\displaystyle{D=C(\widetilde{N})}} & & & &
\ar@{-}[ddddd]^-{\displaystyle{z^{-1}D=C(z^{-1}\widetilde{N})}} \\
& & & & \ar[llll]_-{\displaystyle{h_D}}   \\
& & & & \ar[dll]_-{\displaystyle{h_F}}  \\
& & F=C^{MS}(M_N,f_N,v_N) \ar[llu]_-{\displaystyle{c}}  
\ar[dl]_-{\displaystyle{d_F}} & & \\
& & & & \\
\ar@{-}[rrrr] & & & & }$$
The philosophy here is that ${\mathcal C}(\phi)$ counts the
$\overline{v}$-gradient flow lines of $\overline{f}:\overline{M} \to
\R$, as follows~:
\begin{itemize}
\item[(i)] the $(z^{-1}p,q)$-coefficient of $h_D:z^{-1}D_i \to D_i$ counts the 
number of portions in $M_N$ of the $\overline{v}$-gradient flow lines
which start in $z^{-1}M_N$, enter $M_N$ at $z^{-1}p \in z^{-1}N$,
exit at $q \in N$ and end in $zM_N$,
\item[(ii)] the $(z^{-1}p,q)$-coefficient of $h_F:z^{-1}D_i \to F_i$ counts the 
number of portions in $M_N$ of the $\overline{v}$-gradient flow lines 
which start in $z^{-1}M_N$, enter $M_N$ at $z^{-1}p \in z^{-1}N$ and end at $q \in M_N$,
\item[(iii)] the $(p,q)$-coefficient of $c:F_i \to D_{i-1}$ counts the 
number of portions in $M_N$ of the $\overline{v}$-gradient flow lines 
which start at $p \in M_N$, exit at $q \in N$, and end in $zM_N$.
\end{itemize}
Then for $j=1,2,3,\dots$ 
the $(p,z^jq)$-coefficient of $h_F(h_D)^{j-1}c:F_i \to z^jF_i$ is
the number of the $\overline{v}$-gradient flow lines which start at
$p \in M_N$ and end at $z^jq \in z^jM_N$, crossing the walls 
$N,zN,\dots,z^{j-1}N$. If such is the case, i.e. if the chain map
$h$ is {\it gradient-like} in the terminology of Ranicki \cite{Ranicki1999},
this is just the $(p,z^jq)$-coefficient
of $d_{C^{Nov}(M,f,v)}$, so $\widehat{F}=C^{Nov}(M,f,\alpha)$.
Pajitnov \cite{Pajitnov1999} constructed a $C^0$-dense subspace
$\GCCT(f) \subset \GT(f)$ of gradient-like vector fields $v$ for which
there exist a $CW$ structure $N$ and a gradient-like chain map $h$.
Cornea and Ranicki \cite{CorneaRanicki} construct for any
$v \in GT(f)$ a Morse function $f':M \to S^1$ arbitrarily close to $f$
with $v' \in \GT(f')$ such that
$$C^{Nov}(M,f',v')~=~{\mathcal C}(\phi)~.$$
\indent The projection
$$p~:~C(M;\widehat{\Z[\Pi]})~=~{\mathcal C}(\phi) \to {\rm coker}(\phi)~\cong~
\widehat{F}$$
is a chain equivalence of based f.g. free $\widehat{\Z[\Pi]}$-module
chain complexes, with torsion 
$$\tau(p)~=~\sum\limits^{\infty}_{i=0}(-)^i \tau\big(1-zh_D:
(D_i)_{\alpha}((z)) \to (D_i)_{\alpha}((z))\big) \in K_1(\widehat{\Z[\Pi]})~.$$
If $h$ is a gradient-like chain map the torsion of $p$ is a measure
of the number of closed orbits of the $v$-gradient flow in $M$,
i.e. the closed flow lines $\gamma:S^1 \to M$
(Hutchings and Lee \cite{HutchingsLee1},\cite{HutchingsLee2},
Pajitnov \cite{Pajitnov1999},\cite{Pajitnov2001}, Sch\"utz 
\cite{Schuetz2000},\cite{Schuetz2001}).

The algebraic surgery treatment of high-dimensional
knot theory in Ranicki \cite{Ranicki1998} gives the following algebraic
model for the circle valued Morse function on a knot complement.

\noindent{\it Example.} Let $k:S^n \subset S^{n+2}$ be a knot with
$\pi_1(S^{n+2}\backslash k(S^n))=\Z$.
The complement of a tubular neighbourhood $k(S^n) \times D^2 \subset S^{n+2}$
is an $(n+2)$-dimensional manifold with boundary
$$(M,\partial M)~=~({\rm cl.}\big(S^{n+2} \backslash (k(S^n) \times D^2)\big),
k(S^n) \times S^1)$$
with 
$$\pi_1(M)~=~\Z~,~\pi_1(\overline{M})~=~\{1\}~,~H_*(M)~=~H_*(S^1)~.$$ 
Let $f:(M,\partial M) \to S^1$ be a map representing $1 \in H^1(M)=\Z$,
with $f\vert:\partial M \to S^1$ the projection. Making $f$ transverse
regular at $0 \in S^1$ there is obtained a Seifert surface 
$N^{n+1}=f^{-1}(0) \subset M$ for $k$, with $\partial N=k(S^n)$.
As before, cut $M$ along $N$ to obtain a fundamental domain
$(M_N;N,z^{-1}N)$ for the infinite cyclic cover $\overline{M}=f^*\R$
of $M$. For any $CW$ structures on $N,M_N$ write the reduced chain complexes as
$$\dot C(N)~=~C(N,\{{\rm pt.}\})~~,~~\dot C(M_N)~=~C(M_N,\{{\rm pt.}\})~.$$
The inclusions $G:N \to M_N$, $H:z^{-1}N \to M_N$ induce $\Z$-module  
chain maps 
$$G~:~\dot C(N) \to \dot C(M_N)~,~H~:~z^{-1}\dot C(N) \to \dot C(M_N)$$ 
such that $G-H:\dot C(N) \to \dot C(M_N)$ is a chain equivalence.
The chain map 
$$e~=~(G-H)^{-1}G~:~\dot C(N) \to \dot C(N)$$ 
is a generalization of the Seifert matrix, such that up
to $\Z$-module chain homotopy 
$$1-e~ =~ -(G-H)^{-1}H~:~\dot C(N) \to \dot C(N)$$ 
and such that there is defined a $\Z[z,z^{-1}]$-module chain equivalence
$$C(\overline{M},\R)~\simeq~
C\big(e+z(1-e):\dot C(N)[z,z^{-1}] \to \dot C(N)[z,z^{-1}]\big)~.$$
Let $\dot N={\rm cl.}(N \backslash D^{n+1})$, for any embedding $D^{n+1}\subset
N\backslash \partial N$. For any handlebody decomposition of the
$(n+1)$-dimensional cobordism $(\dot N;k(S^n),S^n)$ with $c_i(N)$ $i$-handles
$$\dot N~=~k(S^n) \times [0,1] \cup 
\bigcup^n_{i=1}\bigcup_{c_i(N)}D^i \times D^{n+1-i}$$
there exists a Morse function $f:M \to S^1$ in the homotopy class
$1 \in [M,S^1]=H^1(M)=\Z$ with
$$c_i(f)~=~c_i(N)+c_{i-1}(N)$$
critical points of index $i$.
In this case the algebraic model for $C^{Nov}(M,f,v)$ has
$$\begin{array}{l}
D~=~C(N)~=~\Z \oplus \dot D~,~\dot D_i~=~\Z^{c_i(N)}~,\\[1ex]
F~=~C^{MS}(M_N,f_N,v_N)~=~{\mathcal C}(e:\dot D \to \dot D)~,\\[1ex]
d_F~=~\begin{pmatrix} d_{\dot D} & e \\ 0 & -d_{\dot D} \end{pmatrix}~:~
F_i~=~\dot D_i \oplus \dot D_{i-1} \to F_{i-1}~=~\dot D_{i-1} \oplus
\dot D_{i-2}~,\\[2ex]
c~=~\begin{pmatrix} 0 & 1 \end{pmatrix}~:~F_i~=~\dot D_i \oplus \dot
D_{i-1} \to D_{i-1}~,\\[1ex]
h_D~=~0~:~z^{-1}D_i \to D_i~,\\[1ex]
h_F~=~\begin{pmatrix} 1-e \\ 0 \end{pmatrix}~:~
z^{-1}D_i \to F_i~=~\dot D_i\oplus \dot D_{i-1}
\end{array}$$
with algebraic Novikov complex
$$\begin{array}{ll}
d_{\widehat{F}}&=~d_F+\sum\limits^{\infty}_{j=1}z^jh_F(h_D)^{j-1}c\\
&=~\begin{pmatrix} d_{\dot D} & e+z(1-e) \\ 0 & -d_{\dot D} \end{pmatrix}~:~
\widehat{F}_i~=~(\dot D_i \oplus \dot D_{i-1})((z)) \to 
\widehat{F}_{i-1}~=~(\dot D_{i-1} \oplus \dot D_{i-2})((z))
\end{array}$$
such that $H_*(\widehat{F})=H^{Nov}_*(M,f)$.
The short exact sequences of $\Z((z))$-modules
$$0 \to H_i(N)((z)) \xrightarrow[]{e+z(1-e)} H_i(N)((z)) \to H^{Nov}_i(M,f) \to 0$$
can be used to express the Novikov numbers $b^{Nov}_i(M,f),q^{Nov}_i(M,f)$
of the knot complement in terms of the Alexander polynomials
$$\Delta_i(z)~=~{\rm det}(e+z(1-e):H_i(N)[z,z^{-1}] \to H_i(N)[z,z^{-1}])
\in \Z[z,z^{-1}]~(1 \leqslant i \leqslant n)~,$$
generalizing the case $n=1$ due to Lazarev \cite{Lazarev}.
For $n \geqslant 4$ and $\pi_1(M)=\Z$ the following conditions are equivalent~:
\begin{itemize}
\item[(i)] the knot fibres, i.e.  $f:M \to S^1$ is homotopic to the
projection of a fibre bundle, with no critical points,
\item[(ii)] $b^{Nov}_*(M,f)=q^{Nov}_*(M,f)=0$, i.e. $H^{Nov}_*(M,f)=0$,
\item[(iii)] the constant and leading coefficients of $\Delta_*(z)$ are 
$\pm 1 \in \Z$.\hfill\qed
\end{itemize}

There is also a more refined version of the algebraic model for circle
valued Morse theory, using the noncommutative Cohn localization
$\Sigma^{-1}A_{\alpha}[z,z^{-1}]$ of $A_{\alpha}[z,z^{-1}]$ inverting
the set $\Sigma$ of square matrices of the form $1-zh$ for a square
matrix $h$ over $A$. Indeed, the formula for the differentials in the 
algebraic Novikov complex 
$$d_{\widehat{F}}~=~d_F+zh_F(1-zh_D)^{-1}c$$
is already defined in $\Sigma^{-1}A_{\alpha}[z,z^{-1}]$.
See Farber and Ranicki \cite{FarberRanicki}
and Ranicki \cite{Ranicki1999} for further details of the
construction. Farber \cite{Farber1999} applied the refinement to
obtain improvements of the Morse-Novikov inequalities, using homology 
with coefficients in flat line bundles instead of Novikov homology.
It should be noted that the natural morphism
$\Sigma^{-1}A_{\alpha}[z,z^{-1}] \to A_{\alpha}((z))$ is 
injective for commutative $A$ with $\alpha=1$, but it is 
not injective in general (Sheiham \cite{Sheiham2001}).


\newpage
\addcontentsline{toc}{section}{References}
\providecommand{\bysame}{\leavevmode\hbox to3em{\hrulefill}\thinspace}

\end{document}


\endRefs
\enddocument